\documentclass[12pt]{article}
\usepackage{defs1}
  \usepackage{slashed}
\title{The universal $\eta$-invariant for manifolds with boundary}
\author{Ulrich Bunke\thanks{NWF I - Mathematik,
Universit{\"a}t Regensburg,
93040 Regensburg,
GERMANY, ulrich.bunke@mathematik.uni-regensburg.de}  
}

\theoremstyle{definition}

\newcommand{\btmr}{\mathbf{TM}}

\newcommand{\Moore}{\mathrm{Moore}}
\newcommand{\asss}{{\tt ass}}
\newcommand{\PD}{{\tt PD}}
\newcommand{\bP}{\mathbf{P}}
\newcommand{\fibre}{\mathrm{fibre}}
\newcommand{\bku}{{\mathbf{ku}}}
\newcommand{\MF}{{\mathcal{MF}}}

\newcommand{\Iso}{\mathrm{Iso}}

\newcommand{\Gr}{{\mathrm{Gr}}}

\newcommand{\bK}{{\mathbf{K}}}

\newcommand{\spec}{{\tt spec}}

\newcommand{\cD}{{\mathcal{D}}}

 \newcommand{\Cone}{{\tt Cone}}
 
 \newcommand{\Vect}{{\tt Vect}}

\renewcommand{\Dirac}{\slashed{D}}

\begin{document}
\maketitle
\begin{abstract}
 We extend the theory of the universal $\eta$-invariant to the case of  bordism groups of manifolds with boundaries. This allows the construction of secondary descendants of the universal $\eta$-invariant.  We obtain an interpretation of Laures' $f$-invariant as an example of this general construction. \textcolor{black}{As an aside  we improve a recent result by Han-Zhang about the modularity of a
certain formal power series of  $\eta$-invariants.}  
\end{abstract}

\setcounter{tocdepth}{1}      
\tableofcontents

\section{Introduction}

In this paper we investigate the question how elements in bordism groups of manifolds with boundary can be detected using spectral invariants of Dirac operators, namely the $\eta$-invariant of Atiyah-Patodi-Singer \cite{MR0397797}. The corresponding problem for bordism groups of closed manifolds has been thouroughly studied in \cite{2011arXiv1103.4217B} and led to the introduction of the universal $\eta$-invariant. The purpose of the present paper is to extend this theory from closed manifolds to manifolds with boundary. One motivation for this generalization is  to prove
Theorem \ref{thm1000} which refines a recent result in \cite[Thm 1.1]{2013arXiv1312.7494H}.

\bigskip

The $\eta$-invariant of a Dirac operator $\Dirac_{M}$ on a closed  manifold $M$ was defined in \cite{MR0397797} as its $\zeta$-regularized signature 
\begin{equation}\label{jdkkwedewd}\eta(\Dirac_{M}):=\sum_{\lambda\in \spec(\Dirac_{M})\setminus \{0\}} \mathrm{mult}(\lambda) \frac{\sign(\lambda)}{|\lambda|^{s}}\:\:\Big|_{s=0}\ .\end{equation}
The sum converges if $\Ree(s)$ is large, and the value at $s=0$ is defined as the evaluation of  the  meromorphic continuation of the sum which happens to be regular at this point. The $\eta$-invariant $\eta(\Dirac_{M})$ is one term in the APS index formula \cite{MR0397797} for  the index of a Dirac operator $\Dirac_{W}$  extending $\Dirac_{M}$ over a zero bordism $W$ of $M$ with   APS boundary conditions:
\begin{equation}\label{gdhwdgwdwd}
\ind(\Dirac_{W})_{APS}=\int_{W} \mbox{index density} -\frac{\eta(\Dirac_{M})+\dim(\ker(\Dirac_{M}))}{2}\ .\end{equation} 
This formula is the starting point for the construction of bordism invariants of $M$. These topological invariants are derived not just from a single $\eta$-invariant but from the relation between the $\eta$-invariants of  the  twisted Dirac operators $\eta(\Dirac_{M}\otimes \bV)$ for various geometric vector bundles $\bV$ on $M$. A geometric vector bundle  is a  triple $\bV=(V,h^{V},\nabla^{V})$ consisting of a complex vector bundle with hermitean metric and metric connection.  
 The rough idea is to form suitable linear combinations of $\eta$-invariants such that the integral of the index density (which encodes the continuous dependence of the $\eta$-invariant on geometric data) drops out, and to consider equivalence classes of the values in $\R/\Z$ in order to get rid of the contribution of the index. A typical example of this idea is the well-known  construction of the $\rho$-invariant $$[\eta(\Dirac_{M}\otimes \bV)-\dim(\bV)\eta(\Dirac_{M})]\in \R/\Z$$ for a flat geometric  bundle $\bV$.
The universal invariant of this kind is the analytic version of the universal $\eta$-invariant $\eta^{an}$ introduced in \cite{2011arXiv1103.4217B}. 
\textcolor{black}{It is called universal since it is indeed the universal case for a certain kind of construction of bordism invariants derived from $\eta$-invariants from which others can be derived by specialization. This was demonstrated in a number of examples in  \cite{2011arXiv1103.4217B}.}
One of the main results in \cite{2011arXiv1103.4217B} is a complete description of $\eta^{an}$ in terms of homotopy theory. To this end we introduce a topological version $\eta^{top}$ of the universal $\eta$-invariant defined in terms Thom spectra and their
$\Q$ and $\Q/\Z$-versions, and we show a secondary index theorem stating that $\eta^{an}=\eta^{top}$.

\bigskip

There are various ways to extend the definition of $\eta(\Dirac_{M})$ to manifolds with boundaries.
In order to extend \eqref{jdkkwedewd} one must define a selfadjoint extension of $\Dirac_{M}$
by choosing suitable boundary conditions. Another possibility   would be to attach an infinite cylinder. After attaching the cylinder the operator $\Dirac_{M}$ has a natural selfadjoint extension. If $M$ is compact one can interpret the sum \eqref{jdkkwedewd}
as a trace of a function of $\Dirac_{M}$.  In the case with boundaries completed by cylinders
the corresponding trace does not exist on the nose but can be defined using a  regularization procedure, e.g. by employing the $b$-calculus of Melrose or \cite{MR2191484}. Finally, one can avoid non-compact manifolds or boundaries at all by
forming doubles. In the present paper we prefer this last method whose details are given in Section \ref{dwdwd1}. 

\bigskip
 
The natural domain for the generalization of the universal $\eta$-invariant $\eta^{an}$ to manifolds with boundary    is a relative bordism group. It can be defined as a homotopy group of a relative bordism spectrum. \textcolor{black}{The data for such a spectrum 
consists of a diagram of spaces \begin{equation}\label{bnxqklxnlqxqwx}\xymatrix{A\ar[rr]^{i}\ar[dr]&&B\ar[dl]\\&BSpin^{c}&}\ .\end{equation}}The relative bordism spectrum $M(B,A)$ will be defined in Section \ref{jqwkdjwqndwkqjd}. The elements of $\pi_{n}
(M(B,A))$ can be interpreted as bordism classes of $n$-dimensional $B$-manifolds with boundary  on which  the $B$-structure is refined to an $A$-structure.  

\bigskip

The topological version of the universal $\eta$-invariant $\eta^{top}$ was  defined for every spectrum in  \cite{2011arXiv1103.4217B} (see Section \ref{fiewefwef}). In particular it can be applied to $M(B,A)$. We know that it can detect  torsion
elements in $\pi_{n}
(M(B,A))$ which survive $K$-localization.

 In the present paper (Sections \ref{jksaxsaxsx}, \ref{dwdwd1}, and  \ref{jdkqdqd}) we generalize the construction of the analytic version $\eta^{an}$  to $M(B,A)$. 
Extending the main result of  \cite{2011arXiv1103.4217B} we  show the secondary index Theorem \ref{duhdiuqwdwqd} stating that $\eta^{top}=\eta^{an}$.

\bigskip

The simplest definition of $\eta^{an}$  given in Definition \ref{lodewdwed} does not inolve $\eta$-invariants on manifolds with boundary but $\R/\Z$-valued indizes on zero bordisms. In order to give a formula in terms of $\eta$-invariants (called the intrinsic formula) we need to choose a further structure called geometrization  \cite{2011arXiv1103.4217B}. In a certain sense the notion of a geometrization generalizes the notion of a connection on a principal bundle. 
Following the lines of  \cite{2011arXiv1103.4217B} in Section \ref{sjkfsfsfsrfsrf} we  extend the notion of a geometrization to the relative case. In    Theorem \ref{widoqdwqd} we provide the corresponding intrinsic formulas for $\eta^{an}$.

The notion of a geometrization involves differential $K$-theory and the construction of differential $K$-theory classes from geometric vector bundles. This construction is  usually called the cycle map. For manifolds with boundary we must introduce the relative version of  differential $K$-theory for pairs $(M,N)$ of a manifold $M$ and a submanifold $N$. 
As a technical ingredient of independent interest  we construct the cycle map for relative differential $K$-theory in Section \ref{klasxasxasx}. It associates a relative differential $K$-theory class to a pair $(\bV,\rho)$ of a geometric bundle $\bV$ on $M$ together with a trivialization $\rho$ of geometric bundles of the restriction $\bV_{|N}$.

\bigskip

In a certain sense the  universal $\eta$-invariant is   a secondary invariant for the index of Dirac operators. The novelty of the case with boundaries is that it allows to define tertiary descendants of the index   which are secondary for the universal $\eta$-invariant (see Section \ref{ascsc}). 
A first example of such an invariant has been studied in \cite{MR2652438} and was identified with Laures' $f$-invariant \cite{MR1660325}, \cite{MR1781277}. This invariant can detect certain elements in the stable homotopy groups of spheres. \textcolor{black}{As an illustration,  in  Section \ref{dsvdsvdsvvvv} of 
 the present paper}  we show how this fits into the general framework of the universal $\eta$-invariant and   how general properties of the universal $\eta$-invariant imply (already known) features of the $f$-invariant. 

\bigskip

Motivated by recent work \cite{MR2153079} in Section \ref{dioqd} we discuss a second example, a $Spin$-bordism version $f^{Spin}$ of the $f$-invariant. In this case, again as an illustration, we specialize our general theory and provide a geometrization and an intrinsic formula for $f^{Spin}$.

 We try to keep this paper short and refer to \cite{2011arXiv1103.4217B} for many details of the language and some arguments.
 
 \bigskip

\textit{Acknowledement: I thank Fei Han and Weiping Zhang for the interesting discussion on their recent paper \cite{2013arXiv1312.7494H}.
One of the purposes of the present  paper, in particular  of Section  \ref{dioqd}, is to answer some questions  asked in \cite{2013arXiv1312.7494H}. \textcolor{black}{I further thank the referee for pointing out various  incorrect statements  in the first version of this paper}.}

\section{The topological universal $\eta$-invariant}\label{fiewefwef}

In this section we introduce the topological version of the universal $\eta$-invariant.
Let $E$ denote  a spectrum, $n$ be an integer, and $\pi_{n}(E)_{tors}\subseteq \pi_{n}(E)$ denote  the torsion subgroup of the homotopy group of $E$ in degree $n$.
 The universal $\eta$-invariant introduced in \cite{2011arXiv1103.4217B}
 is a homomorphism of abelian groups
\begin{equation}\label{xajksx}\eta^{top}:\pi_{n}(E)_{tors}\to Q^{\R}_{n}(E)\ .\end{equation}
   In the following we first describe the target group $Q^{\R}_{n}(E)$ and then the construction of $\eta^{top}$. 
   
   \bigskip
 
 For an abelian group $G$ and a spectrum $E$ we let $\Moore(G)$ denote the Moore spectrum of $G$ (see \cite[Sec. 2]{MR551009}) and abbreviate $EG:=E\wedge \Moore(G)$.  At various places we will  use the fibre sequence of spectra
\begin{equation}\label{flwfwefii}\Sigma^{-1} E\R \to \Sigma^{-1}E\R/\Z \to E \to E\R\to E\R/\Z\ .\end{equation} Let $K$ denote the complex $K$-theory spectrum.
  \begin{ddd}\label{wefuiwefwpefepfeellee}We  define the abelian group
$$Q^{\R}_{n}(E):=\frac{\Hom^{cont}(K^{0}(E),\pi_{n+1}(K\R/\Z))}{U_{n}^{\R}}\ ,$$
where the subgroup   $U_{n}^{\R}\subseteq \Hom^{cont}(K^{0}(E),\pi_{n+1}(K\R/\Z))$ consists of all homomorphisms determined by  elements $y\in \pi_{n+1}(E\R)$ as compositions
 $$ K^{0}(E)\ni \phi \mapsto \left(\Sigma^{n+1}S \stackrel{y}{\to} E\R\stackrel{\phi }{\to} K\R\to K\R/\Z\right)\in \pi_{n+1}(K\R/\Z )\ .$$ 
 \end{ddd}
  In order to talk about continuous homomorphisms we equip  the group $K^{0}(E)$ with the profinite topology \cite[Def. 4.9]{board},  and the group $\pi_{n+1}(\R/\Z)$ with the discrete topology. 
 
 \bigskip

We now construct the homomorphism $\eta^{top}$.
We will describe the value  $\eta^{top}(x)\in Q^{\R}_{n}(E)$ for $x\in \pi_{n}(E)_{tors}$. The fibre sequence \eqref{flwfwefii} induces  a long  exact sequence of abelian groups
$$\pi_{n+1}(E\R)\to \pi_{n+1}(E\R/\Z)\to \pi_{n}(E)\to \pi_{n}(E\R)\ .$$ Since it is torsion the image of $x$ in $\pi_{n}(E\R)$ vanishes.
 Therefore we can find a lift $\tilde x\in \pi_{n+1}(E\R/\Z)$ of $x$ which is unique up to elements coming from $\pi_{n+1}(E\R)$.
The element $\tilde x$ determines a continuous homomorphism
\begin{equation}\label{jdjkededkejjjhw}K^{0}(E)\ni \phi \mapsto \left(\Sigma^{n+1}S \stackrel{\tilde x}{\to} E\R/\Z\stackrel{\phi }{\to}  K\R/\Z\right)\in \pi_{n+1}(K\R/\Z)\ .\end{equation}
\begin{ddd}\label{duedeiduwhedhwed}We define the topological universal $\eta$-invariant    as the map \eqref{xajksx} such that  
 $\eta^{top}(x)$ is represented by the composition \eqref{jdjkededkejjjhw}.
 \end{ddd}
The latter is well-defined independently of the choice of $\tilde x$  
exactly since we take the quotient by $U_{n}^{\R}$ in the definition of $Q^{\R}_{n}(E)$. 

\bigskip

We refer to \cite[Sec. 2.4]{2011arXiv1103.4217B} for an analysis  of $\eta^{top}$ in terms of stable homotopy theory.


\section{Cycles for  relative bordism theory}\label{jqwkdjwqndwkqjd}

A morphism of spaces $B\to  BSpin^{c}$ gives rise to a Thom spectrum $MB$.
The associated generalized homology theory is the bordism theory of $B$-manifolds.
 A \textcolor{black}{diagram \eqref{bnxqklxnlqxqwx}} induces a map of Thom spectra which we will extend to a fibre sequence
\begin{equation}\label{jkklwdqwd}\Sigma^{-1}M(B,A)\to MA\to MB\to  M(B,A)\end{equation}
 in order to define the spectrum
$M(B,A)$ as the cofibre.  
By the Thom-Pontrjagin construction the
 generalized  cohomology theory represented by the spectrum $M(B,A)$ can be described as the  bordism theory of $B$-manifolds with boundary, on which the $B$-structure is refined to an $A$-structure.   
 
 \begin{rem}{\rm
Before we can talk about geometric cycles we must fix once and for all models for the spaces $BSpin^{c}(k)$ and universal bundles $\xi^{k}$, for the maps $s:BSpin^{c}(k)\to BSpin^{c}(k+1)$, and for the isomorphisms $s^{*}\xi^{k+1}\cong \underline{\R}_{BSpin^{c}(k)}\oplus \xi^{k}$. \textcolor{black}{
We can and will assume, after some homotopy equivalent replacement, that the maps  to $BSpin^{c}$ in \eqref{bnxqklxnlqxqwx} are fibrations, and that $i$ is a cofibration. Ay this assumption the map $A\to B$ is the inclusion of a subspace.}

In the  definition of cycles and relations below  we use the Riemannian metric in order to define normal vector fields at the boundary which in turn are required to define the boundary restriction of a $B$-structure. 
Later we will also use the Levi-Civita connection in order to define geometric differential operators.
} \hB\end{rem}


In the following we give a more detailed description of the cycles  $(M,N,f)$ for classes in $\pi_{n}(M(B,A))$ and the equivalence relation.

\begin{enumerate}
\item $M$ is a compact $n$-dimensional Riemannian manifold
with boundary $N$. We assume that the Riemannian metric   has a product structure near the boundary of $M$. 
\item  $f:M\to B$ is a continuous  map which is refined by a stable normal $B$-structure \cite[Def. 3.1]{2011arXiv1103.4217B}. By definition such a refinement consists of a lift $\hat f$ and a homotopy filling of the diagram
\begin{equation}\label{hhdwudiwdwdw}\xymatrix{&&BSpin^{c}(k)\ar[d]\\M\ar[r]_{f}\ar[rru]^{\hat f}&B\ar[r]&BSpin^{c}}\end{equation} together with
  the choice of an isomorphism of real vector bundles
\begin{equation}\label{udiuqdhqwdwqd}TM\oplus \hat f^{*}\xi^{k}\cong \underline{\R^{n+k}}_{M}\ .\end{equation}
The right-hand side  is the notation for the trivial $n+k$-dimensional real vector bundle,
and $\xi^{k}$ is the universal $k$-dimensional real vector bundle with a $Spin^{c}$-structure on $BSpin^{c}(k)$.
\item The restriction $f_{|N}$ has values in the subspace $A$. 
Note that the stable normal $A$-structure on $N$ is given by the restriction $\hat f_{|N}$ and the isomorphism
\begin{equation}\label{djckcecew}TN\oplus (s\circ \hat f_{|N})^{*}\xi^{k+1}\cong TN\oplus \underline{\R}_{N}\oplus \hat f_{|N}^{*}\xi^{k}\cong TM_{|N}\oplus \hat f_{|N}^{*}\xi^{k}\cong  \underline{\R^{n+k}}_{N}\ ,\end{equation}
where $s:BSpin^{c}(k)\to BSpin^{c}(k+1)$ denotes the stabilization, the second isomorphism uses the outer normal vector field at $N$, and  the last isomorphism is given by restriction of 
\eqref{udiuqdhqwdwqd} to the boundary.
\end{enumerate}


A zero bordism  $(W,F)$ of such a cycle is given by the following data:
\begin{enumerate}
\item
A compact $n+1$-dimensional Riemannian manifold $W$ with corners of codimension $2$ and a partition of the boundary $\partial W=\partial_{A} W\cup \partial_{B} W$ such  that
$\partial_{A} W\cap \partial_{B} W$ is the codimension two stratum. We  assume a product structure of the Riemannian metric near all boundary faces and that the faces meet at a right angle in the corner. 
\item  $F:W\to B$  is a map  which is refined to a normal $B$-structure.
\item The restriction of $F_{|\partial_{A}W}$ has values in $A$.
\item There is an isomorphism of cycles
$$\partial_{B}(W,F)\cong (M,N,f)\ .$$
In detail this means that $M\cong \partial_{B}W$ equipped with the  induced Riemannian metric,  
$F_{|N}=f$, and the refinement of $f$ to a normal $B$-structure is induced from that of $F$ similarly as in \eqref{djckcecew}.
 \end{enumerate}

\section{The group $Q^{\R}_{n}(M(B,A))$}\label{jksaxsaxsx}

In Definition \ref{wefuiwefwpefepfeellee} we  introduced the group $Q_{n}^{\R}(E)$ for an arbitrary spectrum $E$.
In the case of a relative bordism spectrum $E=M(B,A)$ we can rewrite the definition of $Q^{\R}_{n}(M(B,A))$ in terms of the pair of spaces $(B,A)$. This alternative picture of $Q^{\R}_{n}(M(B,A))$ will be used in the definition of $\eta^{an}$. 

\bigskip

We consider the Eilenberg-MacLane spectrum $HP\R:= H\R[b,b^{-1}]$ with $\deg(b)=-2$. It represents two-periodic real cohomology. The Chern character induces an isomorphism  $\ch:K\R\stackrel{\sim}{\to} HP\R$. 


We now use that Thom spectra defined through maps to $BSpin^{c}$ (as apposed to $BO$)  are $K$-oriented. In particular, we have  a Thom isomorphism
$$\Thom^{K}:K^{0}(B,A)\stackrel{\cong}{\to}K^{0}(M(B,A))\ .$$

Assume now that $n$ is an odd integer. Then we have an isomorphism
$$\Z\cong \pi_{0}(K)\stackrel{b^{-\frac{n+1}{2}}}{\longrightarrow} \pi_{n+1}(K)$$
which induces an isomorphism
$\R/\Z\stackrel{\cong}{\to} \pi_{n+1}(K\R/\Z)$.

 Using  the Thom isomorphism and this identification    we get an identification
 \begin{equation}\label{r4r4r4rfff}Q^{\R}_{n}(M(B,A))\cong \frac{\Hom^{cont}(K^{0}(B,A),\R/\Z)}{\tilde U^{\R}_{n} }\ . \end{equation}
 In order to describe the subgroup $\tilde U_{n}^{\R}$ we rewrite the definition of $U_{n}^{\R}$ using the Riemann-Roch formula. 
 We consider    the projection 
$$ p_{n+1}:\pi_{n+1}(HP\R ) \cong \R\to \R/\Z\ .$$ 
 The elements of $\tilde U^{\R}_{n}$ are homomorphisms determined by elements $y\in H\R_{n+1}(B,A)$ by the formula
 \begin{equation}\label{suiwsws} K^{0}(B,A)\ni \phi\mapsto p_{n+1} ( \langle y,\Td^{-1}\cup \ch(\phi)\rangle )  \ ,\end{equation} where $\Td^{-1}\in HP\R^{0}(B)$ is the universal Todd class pulled back from $BSpin^{c}$.

\section{The doubling construction}\label{dwdwd1}

\subsection{The double of a $B$-manifold}

In this section we describe a doubling construction.  It will be used to get rid of boundary components of type $A$ in order to simplify the analytic arguments later. We consider a zero bordism
$(W,F)$ as in Section \ref{jqwkdjwqndwkqjd}. In this situation we form the double
$$\cD W:=W\cup_{\partial_{A}W} W^{op}\ ,$$
where we use the subsript $op$ in order to indicate that the right copy has the opposite orientation. The double $\cD W$
 is a Riemannian manifold whose boundary is again a double  $$\partial (\cD W)\cong  \cD (\partial_{B}W)\ .$$
 
 \bigskip
 
 The double $\cD W$ has an induced $B$-structure. Its underlying map
 $\cD F:\cD W\to B$ given by $F$ on both copies of $W$. Furthermore, the analog of the isomorphism \eqref{udiuqdhqwdwqd} is given on 
 $W^{op}$    by 
$$TW^{op}\oplus \hat F^{*}\xi^{k}\cong \underline{\R^{n+k}}_{W^{op}}\stackrel{\epsilon}{\to} \underline{\R^{n+k}}_{W^{op}}\ ,$$
where $\epsilon$ flips the $n$th basis vector, and the first isomorphism is \eqref{udiuqdhqwdwqd}  for the normal $B$-structure of the left copy of $W$.
This map glues with the normal $B$-structure on the left copy in view of the fact, that the glueing for the tangent bundle $T\cD W$ is given by the isomorphism
$$TW_{|\partial_{A}W}\cong T\partial_{A}W\oplus \underline{\R}_{\partial_{A}W^{op}}\stackrel{\id\oplus -1}{\cong} T\partial_{A}W^{op}\oplus \underline{\R}_{\partial_{A}W^{op}}
\cong TW^{op}_{|\partial_{A}W^{op}}\  .$$

\subsection{The double of the spinor bundle}

A refinement of the normal $B$-structure of $W$  to a geometric tangential $Spin^{c}$-structure \cite[3.2]{2011arXiv1103.4217B} 
induces a Dirac bundle $S_{W}$. There is a natural construction of a geometric tangential $Spin^{c}$-structure on the double $\cD W$. Since we are only interested in the associated Dirac operators we will describe this double on the level of Dirac bundles, see \cite[Sec. II, Def. 5.2]{MR1031992}.

The opposite Dirac bundle $S_{W^{op}}$ is obtained from $S_{W}$  by replacing the Clifford multiplication $c:TW\times S_{W}\to  S_{W}$ by its negative.
  We obtain the Dirac bundle on $\cD W$ by glueing
 $S_{W}$ on $W$ with $S_{W^{op}}$ on $W^{op}$ using the isomorphism
 $$S_{W|\partial W}\stackrel{c(\nu)}{\to} S_{W^{op}|\partial W^{op}}$$ given by the Clifford multiplication with the outer normal vector field $\nu$.

 \subsection{The double of relative geometric bundles}\label{uhdiuqwdwqdw}

We  now consider a $\Z/2\Z$-graded geometric bundle $\bU$ on $W$ with a product structure near the boundary and  geometry preserving isomorphism
$$\sigma:\bU^{+}_{|\partial_{A}W}\stackrel{\cong}{\to} \bU^{-}_{|\partial_{A}W}\ .$$  In this situation we define the double
$\cD (\bU,\sigma)$, a geometric bundle on $\cD W$, by glueing
$\bU$ on $W$ with $\bU^{-}\oplus \bU^{-}$ on $W^{op}$ using the isomorphism
$$\bU_{|\partial _{A}W}\stackrel{(\sigma,\id)}{\cong} (\bU^{-}\oplus \bU^{-})_{|\partial_{A}W^{op}}\ .$$

\begin{rem}{\rm
Note, that in contrast to the double $\cD W$ and the spinor bundle 
the double $\cD (\bU,\sigma)$ has no reflection symmetry. One should rather think of $\cD(\bU,\sigma)$ as representing a $K$-theory class in $K^{0}(\cD W,W^{op})$ which corresponds to the class of $(U,\sigma)$ in $K^{0}(W,\partial_{A}W)$ under excision.}\hB 
\end{rem}

\section{The analytic universal $\eta$-invariant}\label{jdkqdqd}

In this section we define the analytic version of the   universal $\eta$-invariant   \begin{equation}\label{xajksx1}\eta^{an}:\pi_{n}(M(B,A))_{tors}\to Q^{\R}_{n}(M(B,A)) \end{equation}
 in the relative case. 
 Let $(M,N,f)$ by a cycle representing a relative bordism class  $x\in \pi_{n}(M(B,A))_{tors}$ as explained  in Section \ref{jqwkdjwqndwkqjd}. Then there exists a non-vanishing integer $\ell$ such that $\ell x=0$. Hence there exists a zero bordism $(W,F)$ such that $\partial_{B} (W,F)\cong \ell (M,N,f)$.

\bigskip

\textcolor{black}{Note that $B$ may be non-compact. In order to represent $K$-theory classes by vector bundles
we choose a compact subspace $B_{c}\subseteq B$ which contains the image of $F$ and set $A_{c}:=A\cap B_{c}$.  Let $\phi\in K^{0}(B,A)$.
We can choose a pair  $(V_{c},\rho_{c})$ of a  $\Z/2\Z$-graded  bundle $V_{c}$  
  on $B_{c}$ and an  
isomorphism $\rho:V^{+}_{c|A_{c}}\to V^{-}_{c|A_{c}}$ such that $[\bV_{c},\rho_{c}]=  \phi_{|(B_{c},A_{c})}\in K^{0}(B_{c},A_{c})$.
Note that we have a map of pairs $f:(M,N)\to (B_{c},A_{c})$.
We define $(V,\rho):=f^{*}(V_{c},\rho_{c})$. We can and will adjust the choices such that these objects are smooth. We further choose a geometry $(h^{V},\nabla^{V})$ with a product structure near $N$ such that $\rho$ preserves the metric and connection.
We assume that the geometry has a product structure near $N$.}

\textcolor{black}{
The pair $(U,\sigma):=F^{*}(V_{c},\rho_{c})$ extends $(V,\rho)$ on the $\ell$ copies of $(M,N)$ in the boundary of $W$. We can choose a geometry on $(h^{U},\nabla^{U})$  such that
 $(\bU,\sigma)_{|\partial _{B}W}\cong \ell(\bV,\rho)$. 
}

\bigskip

The twisted Dirac operator $\Dirac_{\cD W}\otimes \cD (\bU,\sigma)$ on the double $\cD W$ with APS-boundary conditions at $\partial \cD W$ is a Fredholm operator $(\Dirac_{\cD W}\otimes \cD (\bU,\sigma))_{APS}$ which is odd with respect to the \textcolor{black}{$\Z/2\Z$ grading of the graded tensor product of the spinor bundle with $\cD(\bU,\sigma)$}.    As usual, we define
its index $ \ind(\Dirac_{\cD W}\otimes \cD (\bU,\sigma))_{APS}\in \Z$ graded dimension of its kernel.  We now consider the quantity
 \begin{equation}\label{jkehehfwiefewff}e:=[\frac{1}{\ell} \ind(\Dirac_{\cD W}\otimes \cD (\bU,\sigma))_{APS}]\in \R/\Z\ .\end{equation}

One checks the following properties in a similar manner as in \cite[Prop 3.4]{2011arXiv1103.4217B}
\begin{enumerate}
\item Using the continuous dependence of $e$ on the geometric data we get the independence of $e$ from  the geometric structures on $M$, $W$, $\bV$ and $\bU$. \textcolor{black}{If we fix the cycle for the relative bordism class, then $e$ only depends on the pair $(V_{c},\rho_{c})$. }
\item We now consider $e$ as a function of the pair  $(V_{c},\rho_{c})$. It is additive under direct sum and depends on the choice of $\rho_{c}$ only up to homotopy. 
Consequently  it only depends on the class $[V_{c},\rho_{c}]=\phi_{|(B_{c},A_{c})}\in K^{0}(B_{c},A_{c})$. We conclude that that
$e$ induces a  homomorphism $\tilde e:K^{0}(B,A)\to \R/\Z$ which is continuous since it factorizes over the restriction along $(B_{c},A_{c})\to (B,A)$.  
\item The class $[\tilde e]\in Q^{\R}_{n}(M(B,A))$  (using   the picture   \eqref{r4r4r4rfff} of $Q^{\R}_{n}(M(B,A))$) of the homomorphism $\tilde e$ does not depend on the choice of the integer $\ell$ and the zero bordism $(W,F)$.
\textcolor{black}{Indeed, if $\tilde e^{\prime}$ is defined for different choices $\ell^{\prime}$ and $(W^{\prime},F^{\prime})$, then we argue as in \cite[proof of Prop. 3.4 (3)]{2011arXiv1103.4217B}.
In a first step we can adjust the choices such that $\ell=\ell^{\prime}$. 
 Using the APS index theorem  we now see that the difference $\tilde e-\tilde e^{\prime}$  
belongs to the group $\tilde U^{\R}_{n}$ appearing in \eqref{r4r4r4rfff}.}
\item The class $[\tilde e]\in Q^{\R}_{n}(M(B,A))$  only depends on the bordism class  $x$.
\textcolor{black}{This is completely analoguous to 
\cite[Prop 3.4, (4)]{2011arXiv1103.4217B}}.

\end{enumerate}

\begin{ddd}\label{lodewdwed}
We define the value of the analytic version of the universal $\eta$-invariant on $x$ by $$\eta^{an}(x):=[\tilde e]\ .$$
\end{ddd}

\section{The index theorem}

Let $\eta^{top}$ be the topological universal $\eta$-invariant defined in Definition \ref{duedeiduwhedhwed}    for $E=M(B,A)$, and $\eta^{an}$ be  the analytical universal $\eta$-invariant defined in Definition  \ref{lodewdwed}.

\begin{theorem}[Secondary index theorem]\label{duhdiuqwdwqd}
We have the equality $\eta^{an}=\eta^{top}$.
\end{theorem}
\proof
We adapt the proof given for the absolute case in \cite{2011arXiv1103.4217B}. 
The remainder of the present section is devoted to the proof of this theorem.

\subsection{A geometric cycle for  $\tilde x$}

We define the pointed space
$C_{\ell}$ as the cofibre of  the $\ell$-fold covering 
\begin{equation}\label{ljdzwudguwzd}S^{1}\stackrel{\ell}{\to} S^{1}\to C_{\ell}\ .\end{equation}
It is a Moore space and related with the Moore spectrum of $\Z/\ell \Z$ by an  equivalence  of spectra $\Sigma^{\infty}C_{\ell}\simeq \Sigma\Moore(\Z/\ell \Z)$. We use the equivalence of spectra
$$M(B,A)\wedge C_{\ell}\simeq M(B\times C_{\ell},A\times C_{\ell}\cup B\times *_{C_{\ell}})$$
in order to interpret elements in the homotopy of $M(B,A)\wedge C_{\ell}$ geometrically as in Section \ref{jqwkdjwqndwkqjd}.

We consider the cofibre sequence of spectra obtained by forming the smash product of the cofibre sequence \eqref{ljdzwudguwzd} with $M(B,A)$. It  induces a long exact sequence in homotopy.
We consider the following segment of this   sequence: $$\pi_{n+2}(M(B,A)\wedge C_{\ell})\stackrel{\partial}{\to} \pi_{n}(M(B,A))\stackrel{\ell}{\to} \pi_{n}(M(B,A))\ .$$ 

Let $x\in \pi_{n}(M(B,A))$ be an $\ell$-torsion element.
Then  we can choose a lift $$\tilde x_{\ell}\in \pi_{n+2}(M(B,A)\wedge C_{\ell})$$  of $x$. It induces a choice of $\tilde x\in \pi_{n+1}(M(B,A)\R/\Z)$ used in the definition of $\eta^{top}(x)$ in Section  \ref{fiewefwef}   via the map $$ M(B,A)\wedge \Sigma^{-2}C_{\ell}\simeq \Sigma^{-1}M(B,A)\Z/\ell \Z\to \Sigma^{-1}M(B,A)\R/\Z\ .$$

\bigskip


Let $(M,N,f)$ and $(W,F)$ be as in Section \ref{jqwkdjwqndwkqjd}.  
We form the Riemannian manifold with boundary $$\tilde W:=(S^{1}\times W)\cup_{S^{1}\times \partial_{B}W\cong \ell  (S^{1}\times M)^{op}} (S_{\ell}^{2}\times M)^{op}\ ,$$
where $S^{2}_{\ell}$ is a two-sphere with $\ell$     \textcolor{black}{open discs  with  pairwise disjoint closures} removed (see \cite[Sec. 3.5]{2011arXiv1103.4217B}). The manifold $S^{2}_{\ell}$ is equipped with a Riemmanain metric with product structure which induces the standard metric on the $\ell$ copies of $S^{1}$ in its boundary. 
We  define a map $\tilde F:\tilde W\to B$ such that its restrictions to the summands are given by
$$S^{1}\times W\stackrel{\pr_{W}}{\to} W\stackrel{F}{\to}B\ , \quad S_{\ell}^{2}\times M\stackrel{\pr_{M}}{\to} M\stackrel{f}{\to} B\ .$$ 
Note that the restriction of $F$ to the boundary of $\tilde W$ factorizes over $A$.
We use the stable framings of $S^{1}$ and $S^{2}_{\ell}$ and the normal $B$-structures on $f$ and $F$ in order to refine the restrictions of $\tilde F$ to the left and right pieces to normal $B$-structures. We refer to \cite[Sec. 3.5]{2011arXiv1103.4217B} for more details. The two refinements can be glued to a normal $B$-structure for $\tilde F$  by a similar construction as for the double in Section \ref{dwdwd1}.  


We furthermore define a map
$\tilde G:\tilde W\to C_{\ell}$ such that its restrictions to the summands are given by 
$$S^{1}\times W \stackrel{\pr_{S^{1}}}{\to} S^{1}\to C_{\ell}\ , \quad S_{\ell}^{2}\times M\stackrel{\pr_{S^{2}_{\ell}}}{\to} S^{2}_{\ell} \stackrel{g}{\to}C_{\ell}\ ,$$
where $g$ is defined as in \cite[(40)]{2011arXiv1103.4217B}.
The geometric cycle $(\tilde W,\partial\tilde W,(\tilde F,\tilde G))$ represents an element $\tilde x_{\ell}\in \pi_{n+2}(M(B,A)\wedge  C_{\ell})$.

\begin{lem}
We have $\partial \tilde x_{\ell}=x$.
\end{lem}
 \proof
 This is shown exactly as in the proof of  \cite[Lemma 3.7]{2011arXiv1103.4217B}.\hB

\subsection{An analytic picture of the pairing
$\langle  \Thom^{K}(\phi), \varepsilon(\tilde x_{\ell})\rangle$.}

Let $\varepsilon:S\to K$ be the unit of the ring spectrum $K$. We get a class
$$\varepsilon(\tilde x_{\ell})\in K_{n+2}(M(B,A)\wedge C_{\ell})\ .$$
For  $\phi\in K^{0}(B,A)$ we consider the pairing
$$\langle \Thom^{K}(\phi), \varepsilon(\tilde x_{\ell})\rangle\in K_{n+2}(C_{\ell})\ .$$ The goal of this subsection is the construction of a geometric representative of this $K$-homology class.
 
 \bigskip 
 
In the following we use the notation $(\bU,\sigma)$ and $(\bV,\rho)$ as in Section \ref{jdkqdqd}. 
We let $\tilde \bV$  be  the $\Z/2\Z$-graded geometric  bundle on $\tilde W$ which is naturally given by
$\pr_{W}^{*}\bU$ on $S^{1}\times W$, and by $\pr_{M}^{*} \bV$ on $S^{2}_{\ell}\times M$. It comes with a natural isomorphism $\tilde \rho:\tilde \bV^{+}_{|\partial \tilde W}\stackrel{\cong}{\to} \tilde \bV^{-}_{|\partial \tilde W}$ induced by the isomorphisms $\sigma$ and $\rho$.

The twisted Dirac operator on 
$\Dirac_{\cD(\tilde W)}\otimes \cD(\tilde \bV,\tilde \rho)$
gives rise to a Kasparov $K$-theory class $$[\Dirac_{\cD(\tilde W)}\otimes \cD(\tilde \bV,\tilde \rho)]\in KK_{n+2}(C(\cD\tilde W),\C)$$
 and thus to a $K$-homology class
$$\cD\tilde G_{*}[\Dirac_{\cD(\tilde W)}\otimes \cD(\tilde \bV,\tilde \rho)]\in KK_{n+2}(C(C_{\ell}),\C)\cong K_{n+2}(C_{\ell}) .$$ 
Here for  pointed space $(X,*)$ we let $C(X)$ denote the algebra of complex-valued continuous functions vanishing at $*$. \textcolor{black}{We further use the representation of the reduced $K$-homology of $X$ in terms of $KK$-theory
$K_{*}(X)\cong KK_{*}(C(X),\C)$. 
\begin{rem}{\rm In the present paper we use the homotopy theoretic definition
$K_{*}(Y):=\pi_{*}(K\wedge Y_{+})$ of the $K$-homology of a space $Y$. The comparison with the $KK$-theoretic version of $K$-homology is accomplished via the geometric $K$-homology  $K^{geom}_{*}(Y)$ defined by \cite{MR679698}. \\
One constructs a natural isomorphism   $K^{geom}_{*}(Y)\stackrel{\sim}{\to} K_{*}(Y)$ using the fact that the geometric cycles carry a fundamental class for the homotopy theoretic version of $K$-homology. First of all, by the Pontrjagin-Thom construction, a geometric cycle carries a homotopy theoretic $Spin^{c}$-bordism theory fundamental class. We obtain the $K$-theory fundamental class by applying the Atiyah-Bott-Shapiro orientation.\\
 Using the Dirac operators on the geometric cycles
one furthermore constructs an isomorphism
$K^{geom}_{*}(Y)\stackrel{\sim}{\to}  KK_{*}(C(Y),\C)$ from the geometric version of $K$-homology to the $KK$-theoretic version. We refer to \cite{MR2330153} for details. 
}
\end{rem}
}

\begin{prop}
We have an equality
$$\cD\tilde G_{*}[\Dirac_{\cD(\tilde W)}\otimes \cD(\tilde \bV,\tilde \rho)]= \langle  \Thom^{K}(\phi), \varepsilon(\tilde x_{\ell})\rangle$$
\end{prop} 
\proof
Under the Thom isomorphism
$$\Thom_{K}:K_{n+2}(M(B,A)\wedge C_{\ell})\cong K_{n+2}(B/A\wedge C_{\ell})$$
the class $\epsilon(x_{\ell})$ corresponds to a class
$$\Thom_{K}(\epsilon(x_{\ell}))\in  K_{n+2}(B/A\wedge C_{\ell})\cong KK_{n+2}(C(B/A\wedge C_{\ell}),\C)\ .$$ We first represent this class in terms of Dirac operators.

In order to define $K$-homology classes associated to Dirac operators on manifolds with boundary we get rid of boundary components by implicitly  completing  the manifolds with infinite cylinders.
The $B$-structure on $\tilde W$ induces a $Spin^{c}$-structure. We choose an extension of the Levi-Civita connection on  $\tilde W$ to a $Spin^{c}$-connection   with a product structure at the boundary.   
The $Spin^{c}$-Dirac operator on $\tilde W$ then gives   rise to a class
$[\Dirac_{\tilde W}]\in KK_{n+2}(C(\tilde W/\partial \tilde W),\C)$. Here
$C(\tilde W/\partial \tilde W)$ is the algebra of continuous functions on $\tilde W$ which vanish on $\partial \tilde W$.
The  map
$$(\tilde F,\tilde G):\tilde W\to  B\times C_{\ell}$$ induces a map
$$(\tilde F,\tilde G)^{*}:C(B/A\wedge C_{\ell})\to C(\tilde W/\partial \tilde W)\ .$$

The element $\Thom_{K}(\varepsilon(\tilde x_{\ell}))$ is given by
$$(\tilde F,\tilde G)_{*}[\Dirac_{\tilde W}] \in KK_{n+2}(C(B/A\wedge C_{\ell}),\C)\ .$$
 
 The argument is similar to that \cite[Lemma 3.8]{2011arXiv1103.4217B} using that $[\Dirac_{\tilde W}]$
is the relative $K$-theory fundamental class of the $Spin^{c}$-manifold  with boundary $(\tilde W,\partial \tilde W)$.


As in \cite[Sec. 3.5]{2011arXiv1103.4217B} one now checks  
 that
\begin{equation}\label{uzgduiqwdqwdqwd}\langle \Thom^{K}(\phi), \varepsilon(\tilde x_{\ell})\rangle =\tilde G_{*}([\Dirac_{\tilde W}] \cap \tilde F^{*}\phi)\in KK_{n+2}(C(C_{\ell}),\C)\ .\end{equation} We have the equality 
$[\tilde \bV,\tilde \rho]=\tilde F^{*}\phi\in K^{0}(\tilde W,\partial \tilde W)$.   Under the isomorpism 
$K^{0}(\tilde W,\partial \tilde W)\cong KK(\C,C(\tilde W/\partial \tilde W))$
the class $[\tilde \bV,\tilde \rho]$ is represented by the Kasparov module
$(C(\tilde W,\tilde \bV), F_{\tilde \rho})$, where
$F_{\tilde \rho}\in \Gamma(\tilde W,\End(\tilde \bV))$ is any extension of $\tilde \rho$ to all of $\tilde W$.
The cap product in \eqref{uzgduiqwdqwdqwd} is represented by
the Kasparov product
$[C(\tilde W,\tilde \bV), F_{\tilde \rho}] \otimes_{C(\tilde W/\partial \tilde W)}[\Dirac_{\tilde W}]$ which can be represented by the Callias type operator \cite{MR1348799}
$$[\Dirac_{\tilde W}\otimes \tilde \bV+ F_{\tilde \rho}]\in KK_{n+2}(C(\tilde W/\partial \tilde W),\C)\ .$$
It is now a consequence of the relative index theorem that
$$\tilde G_{*} [\Dirac_{\tilde W}\otimes \tilde \bV+ F_{\tilde \rho}]=\cD \tilde G_{*}[\Dirac_{\cD \tilde W}\otimes \cD(\tilde \bV,\tilde \rho)]\ .$$ \hB

\subsection{The final step}

In the final step of the proof of Theorem \ref{duhdiuqwdwqd} we must show that
$$[\frac{1}{\ell} \ind(\Dirac_{\cD W}\otimes \cD (\bU,\sigma))_{APS}]\in \R/\Z$$ is equal to the image of
$\cD \tilde G_{*}[\Dirac_{\cD \tilde W}\otimes \cD(\tilde \bV,\tilde \rho)]\in K_{n+2}(C_{\ell})$ under the natural map
$K_{n+2}(C_{\ell})\cong \Z/\ell \Z \to \R/\Z$. 
But this is exactly the fact shown in at the end of the proof of \cite[Theorem 3.6]{2011arXiv1103.4217B}. 
\hB

\section{Relative differential $K$-theory and cycles}\label{klasxasxasx}

 The definition of a geometrization involves differential   $K$-theory, in particular the functor $\hat K^{0}$. More precisely, it employs the Hopkins-Singer version of differential  $K$-theory.  We refer to \cite[Sec. 4.2]{2011arXiv1103.4217B} and the discussion below for a review of the relevant structures. By now there are  various constructions of this version of differential $K$-theory. First of all we have the Hopkins-Singer  construction \cite{MR2192936}. Other, more geometric models are based on   families of Dirac operators   \cite{MR2664467} or structured vector bundles \cite{MR2732065}. All of them give an equivalent functor $\hat K^{0}$ by \cite{MR2608479}.

 In the present paper we need the relative version of differential $K$-theory. Ad-hoc constructions of relative differential cohomology theories have been considered e.g. in \cite{2014arXiv1401.1029F} or \cite{2013arXiv1310.2851B}. But if one represents differential cohomology in terms of sheaves of spectra on the site of smooth manifolds with open covering topology, then the definition of the relative groups becomes completely natural. 
 Therefore we will use this set-up which was developed in detail in \cite{2012arXiv1208.3961B}, see also \cite{2013arXiv1311.3188B}.  
In particular we have a sheaf of spectra $\hat\bK$ (which we will describe in \eqref{jdkwdwdwwwww} below)   representing differential $K$-theory in the sense that $$\hat K^{0}(M)=\pi_{0}(\hat \bK(M))\ .$$ 
 
 \bigskip
 
The evaluation of  the periodic de Rham complex $\Omega P$ on a manifold $M$ is defined  by   $$\Omega P(M):=\Omega(M)[b,b^{-1}]\ ,$$ where $\deg(b)=-2$. 
By $$\sigma^{\ge 0}\Omega P(M)\subset \Omega P(M)$$ we denote its (stupid) truncation  which just neglects the part of negative total degree. Using the Eilenberg-MacLane functor $H$ from chain complexes to spectra we can define a sheaf of spectra $H(\sigma^{\ge 0} \Omega P)$.  The sheaf $\hat \bK$ is now defined as the pull-back
of sheaves of spectra
\begin{equation}\label{jdkwdwdwwwww}\xymatrix{\hat \bK\ar[d]^{I}\ar[r]^{R}&H(\sigma^{\ge 0}\Omega P)\ar[d]\\\underline{K}\ar[r]^{\ch}&\underline{H\R P}}\ .\end{equation} Here the lower horizontal map is the map of constant sheaves of spectra induced by the Chern character $\ch:K\to  HP\R$. 
Furthermore,
the right vertical map is the composition of the map obtained by applying $H$ to the   embedding 
$\sigma^{\ge 0} \Omega P\hookrightarrow\Omega P$ with a version of the de Rham isomorphism $H(\Omega P)\simeq \underline{H\R P}$. We refer to 
\cite{2012arXiv1208.3961B} for the technical details.

\textcolor{black}{
\begin{rem}
{\rm In order to understand this definition  note that in the background we have fixed some $\infty$-category of spectra $\Sp$. The $\infty$-category of sheaves of spectra is the full subcategory of the $\infty$-category $\Fun(\Mf^{op},\Sp)$ of functors from manifolds ($1$-categories are considered as $\infty$-categories using nerves) to spectra which satisfy descent with respect to open coverings. The pull-back \eqref{jdkwdwdwwwww} is understood in this $\infty$-categorial world.  }
\end{rem}}

\bigskip

The short exact sequence of sheaves of complexes $$0\to \sigma^{\ge 0}\Omega P\to \Omega P\to \sigma^{\le-1} \Omega P\to 0$$ induces a fibre sequence of sheaves of spectra
$$\Sigma^{-1} H(\sigma^{\le-1} \Omega P) \stackrel{\partial}{\to} H(\sigma^{\ge 0}\Omega P)\to 
H(\Omega P)\to H(\sigma^{\le-1} \Omega P)\ .$$
We have natural isomorphisms $$\pi_{0}(\Sigma^{-1} H(\sigma^{\le-1} \Omega P)(M)) \cong   \Omega P^{-1}(M)/\im(d)\ ,\quad 
\pi_{0}(H(\sigma^{\ge 0}\Omega P)(M))\cong \Omega P_{cl}^{0}(M)\ ,$$ where $\Omega P_{cl}^{0}(M)\subseteq \Omega P^{0}(M)$ is the subspace of closed forms of total degree zero. Under these isomorphisms the boundary operator $\partial$  induces, after application of $\pi_{0}$,  the de Rham differential $$d:\Omega P^{-1}(M)/\im(d)\to \Omega P_{cl}^{0}(M)\ .$$

 The maps $R$ and $I$ in \eqref{jdkwdwdwwwww} induce, after applying $\pi_{0}$,  the curvature map  and the underlying class map $$R:\hat K^{0}(M) \to \Omega P_{cl}^{0}(M)\ , \quad I:\hat K^{0}(M)\to  K^{0}(M)\ ,$$
where the target of the latter is identified using the natural isomorphism $\pi_{*}(\underline{K}(M))\cong K^{-*}(M)$.  Furthermore, since
\eqref{jdkwdwdwwwww} is cartesian, the fibres of the left and right vertical maps coincide. We thus obtain a fibre sequence of sheaves of spectra
\begin{equation}\label{frfrfrccwcwwsss44}\Sigma^{-1} H(\sigma^{\le-1} \Omega P)\stackrel{a}{\to} \hat \bK\stackrel{I}{\to} \underline{K}\to H(\sigma^{\le-1} \Omega P)\end{equation}
which, after applying $\pi_{0}$ and using    that 
$\pi_{0}(H(\sigma^{\le-1} \Omega P)(M))=0$, gives the exact sequence 
\begin{equation}\label{kdjedklejdlked}K^{-1}(M)\stackrel{\ch}{\to} \Omega P^{-1}(M)/\im(d) \stackrel{a}{\to} \hat K^{0}(M)\stackrel{I}{\to} K^{0}(M)\to 0\ .\end{equation}
 
\bigskip 

We now  generalize these calculations to the  relative case. We consider an embedding of a submanifold
  $i:N\to M$.  Then we define the relative differential $K$-theory group  by
\begin{equation}\label{kwldkwdw}\hat K^{0}(M,N):=\pi_{0}(\fibre(\hat \bK(M)\to \hat \bK(N)))\ .\end{equation}


From the long exact sequence in homotopy we get   a natural isomorphism
$$\pi_{0}\left(\fibre: H(\sigma^{\ge 0}\Omega P)(M)\to H(\sigma^{\ge 0}\Omega P)(N)\right)\cong \Omega P^{0}_{cl}(M,N)\ ,$$
where $\Omega P^{0}_{cl}(M,N)\subseteq \Omega P_{cl}^{0}(M)$ the subspace of all closed forms 
whose restriction to $N$ vanishes. 
 We conclude that in the relative case the curvature becomes a map
 $$R:\hat K^{0}(M,N)\to \Omega P^{0}_{cl}(M,N)\ .$$
In order to generalize the exact sequence \eqref{kdjedklejdlked} to the relative case we
 calculate, term by term, the fibre
 of the evaluation of \eqref{frfrfrccwcwwsss44}
 on the inclusion $N\to  M$. \begin{enumerate}
 \item The homotopy group $\pi_{*}$ of the  fibre of $\underline{K}(M)\to \underline{K}(N)$ is the relative $K$-theory group $K^{-*}(M,N)$. 
 \item $\pi_{0}$ of the fibre of $\hat \bK(M)\to \hat \bK(N)$ is, by definition, $\hat K^{0}(M,N)$.
 \item
 We represent the  fibre of $H(\sigma^{\le-1} \Omega P)(M)\to H(\sigma^{\le-1} \Omega P)(N)$ by
 $$\Sigma^{-1} H\left(\Cone\left( \sigma^{\le-1} \Omega P(M)\to \sigma^{\le-1} \Omega P(N)\right)\right)$$
  Explicitly, this  cone is the complex 
 $$ \sigma^{\le-1} \Omega P^{*}(M)\oplus  \sigma^{\le-1} \Omega P^{*-1}(N)\ , \quad d(\alpha,\beta):= (d\alpha, \alpha_{|N}-d\beta  )\ .$$
 In particular, its cohomology in degree $-1$ is the group $$A^{0}:=\frac{\{(\alpha,\beta)\in \Omega P^{-1} (M)\oplus \Omega P^{-2}(N)\:|\: \alpha_{|N}=d\beta \}}{\{(d\gamma,-d\delta+\gamma_{|N})\:|\: (\gamma,\delta)\in \Omega P^{-2}(M)\oplus \Omega P^{-3}(N)\}}\ .$$
 \end{enumerate}
We denote by  $[\alpha,\beta]$ the class in $A^{0}$ represented by the pair $(\alpha,\beta)$.
The following Lemma is now an immediate consequence of these calculations.
\begin{lem}\label{qdlqwd}
We have an exact sequence
\begin{equation}\label{dwdlwd}K^{-1}(M,N)\stackrel{\ch}{\to} A^{0}\stackrel{a}{\to} \hat K^{0}(M,N)\stackrel{I}{\to}  K^{0}(M,N)\to 0\ .\end{equation}
Furthermore, for $[\alpha,\beta]\in A^{0}$ have
\begin{equation}\label{jhdjkaw}R(a([\alpha,\beta]))=d\alpha\ .\end{equation}
\end{lem}

\textcolor{black}{By abuse of notation (e.g. in the stament of Lemma \ref{dwedewldew}) below we can and will apply $a$ to relative cohomology classes in $H\R P^{-1}(M,N)$.}

\bigskip

 We 
consider a  pair $(\bV,\rho)$ of a $\Z/2\Z$-graded complex vector bundle   over $M$  and a geometry preserving isomorphism
$\rho:\bV_{|N}^{+}\to \bV^{-}_{|N}$. This pair represents a relative $K$-theory class $[\bV,\rho]\in K^{0}(M,N)$. We want to refine this class to a differential $K$-theory class \textcolor{black}{$\widehat{[\bV,\rho]}\in \hat K^{0}(M,N)$. The association
$(\bV,\rho)\mapsto \widehat{[\bV,\rho]}$ is called the cycle map.}

In the following lemma functoriality means that the construction commutes with pull-backs along smooth maps between manifolds.
Note that $\ch(\nabla^{V})\in \Omega P_{cl}^{0}(M,N)$.
\begin{lem}\label{edowedewd}
There exists a functorial construction of  a class
$\widehat{[\bV,\rho]}\in \hat K^{0}(M,N)$ such that
$I(\widehat{[\bV,\rho]})=[\bV,\rho]$ and $R(\widehat{[\bV,\rho]})=\ch(\nabla^{V})$.
\end{lem}
\proof
This was exercise \cite[Ex. 4.180]{2012arXiv1208.3961B}.  Here is the solution.

We use the sheaf $\widehat{\bku}^{\nabla}$ of spectra on smooth manifolds introduced in Section \cite[Sec. 6]{2013arXiv1311.3188B}. 
It  is constructed by group-completing the nerve $\Nerve(\Iso(\Vect^{\nabla}_{\oplus}))$ of the symmetric monoidal stack of vector complex bundles with connections $\Iso(\Vect^{\nabla}_{\oplus})$. It is universal for additive characteristic classes for vector bundles with connection. In particular,
 in \cite[Sec. 6.1]{2013arXiv1311.3188B}  we have constructed a map of sheaves of spectra
$$\hat r:\widehat{\bku}^{\nabla}\to \hat \bK\ .$$  

We start with the construction of the cycle map in the absolute case.
To this end we consider a  geometric bundle $\bV$ on a manifold $M$ as an object
of $\Vect^{\nabla}_{\oplus}(M)$ and therefore as a point in
$\Omega^{\infty} \bku^{\nabla}(M)$. We let
$[\bV]_{\bku^{\nabla}}\in \pi_{0}(\bku^{\nabla}(M))$ be the class of its connected component. Then
$\widehat{[\bV]}:=\hat r([\bV]_{\bku^{\nabla}})\in \hat K^{0}(M)$ is 
the Hopkins-Singer differential $K$-theory class of the bundle $\bV$. The association $\bV\mapsto \widehat{[\bV]}$ is called the cycle map.

\begin{rem}{\rm 
In the models of differential $K$-theory developed in \cite{MR2664467} or \cite{MR2732065} the differential $K$-theory class of a geometric vector bundle is tautologically defined. In the present paper we need the detour  over $\widehat{\bku}^{\nabla}$  since we use a different homotopy theoretic definition of Hopkins-Singer differential $K$-theory in terms of the sheaf $\hat \bK$ which is not immediately related to vector bundles. Recall that the use of sheaves of spectra was essential for the definition of relative differential $K$-theory in \eqref{kwldkwdw}. 
 }\hB
 \end{rem}
 
We now extend the cycle map to the relative case.
We again first construct a class $\widehat{[\bV,\rho]}_{\widehat{\bku}^{\nabla}}\in\pi_{0}( \widehat{\bku}^{\nabla}(M,N))$ and then  set
$$\widehat{[\bV,\rho]}:=\hat r(\widehat{[\bV,\rho]}_{\widehat{\bku}^{\nabla}}).$$

We can consider the isomorphism
$\rho:\bV^{+}_{|N}\to \bV^{-}_{|N}$ as a path  in
$\Nerve(\Iso(\Vect^{\nabla}_{\oplus}(N)))$.  We now apply the group completion map $c:\Nerve(\Iso(\Vect_{\oplus}^{\nabla}(N)))\to \Omega^{\infty} \widehat{\bku}^{\nabla}(N)$ and obtain a path in $\Omega^{\infty} \widehat{\bku}^{\nabla}(N)$. We can consider
the pair
$\left(c(\bV^{+})-c(\bV^{-}),c(\rho)-c(\bV_{|N}^{-})\right)$ of a point and a path
as a point in the standard model of the homotopy fibre $\Omega^{\infty} \bku^{\nabla}(M,N)$ of the restriction map
$\Omega^{\infty} \widehat{\bku}^{\nabla}(M)\to \Omega^{\infty} \widehat{\bku}^{\nabla}(N)$.
By definition, this point represents the class
$\widehat{[\bV,\rho]}_{\widehat{\bku}^{\nabla}}\in \pi_{0}(\widehat{\bku}^{\nabla}(M,N))$. 

\section{Geometrizations}\label{sjkfsfsfsrfsrf}

 The main ingredient of the intrinsic formula for the universal $\eta$-invariant to be discussed in Section \ref{jcknkwjebckewc} is the notion of a geometrization. This new concept was 
  introduced in \cite[Definition 4.3]{2011arXiv1103.4217B}.  In the following we extend the notion of  a geometrization to the relative case.

\bigskip

We consider a   quadruple $(M,N,f,\tilde \nabla^{TM})$, where $(M,N,f)$ is as in  Section \ref{jqwkdjwqndwkqjd}.   and 
$\tilde \nabla^{TM}$ is a $Spin^{c}$-extension of the Levi-Civita connection. It gives rise to the form
$\Td(\tilde \nabla^{TM})\in \Omega P_{cl}^{0}(M)$ representing the class
$f^{*}\Td^{-1}$.

\bigskip

 If $$\cG:K^{0}(B,A)\to \hat K^{0}(M,N)$$ is a continuous \textcolor{black}{homomorphism}, then by the same argument as in the proof \cite[Lemma 4.2]{2011arXiv1103.4217B}
  there exists a continuous \textcolor{black}{homomorphism} $ c_{\cG}$,   called \textcolor{black}{a} cohomological character, which completes the following diagram: $$\xymatrix{&\hat K^{0}(M,N)\ar[dr]^{\Td(\tilde \nabla^{TM})\wedge R(-)}&\\K^{0}(B,A)\ar[ru]^{\cG}\ar[rd]^{\Td^{-1}\cup \ch}&&\Omega P_{cl}^{0}(M,N)\\&H\R P^{0}(B,A)\ar@{..>}[ru]^{c_{\cG}}&}\ .$$
  
The periodic de Rham complex $\Omega P$ and the periodic real cohomology $H\R P$ have an additional grading which counts the power of the variable $b$.
 We will call this the $b$-degree. For example, if $\omega\in \Omega^{5}(M)$, then $b^{-3}\omega\in \Omega P(M)$ has cohomological degree $11$ and $b$-degree $-3$. 
 \begin{ddd}\label{iuolqdqdwd}
A geometrization of $(M,N,f,\tilde \nabla^{TM})$ is a continuous \textcolor{black}{homomorphism} $\cG:K^{0}(B,A)\to  \hat K^{0}(M,N)$ \textcolor{black}{which admits a $b$-degree-preserving cohomological character.}
\end{ddd} 

The construction of geometrizations is a non-trivial matter. Here we demonstrate such a construction 
 in  the example     where
  $A=*$, $B= BSpin$, and where $B\to BSpin^{c}$
is the canonical map. This example will be employed in Section \ref{dioqd}.
Thus we consider  a Riemannian manifold  $M$ with an embedded submanifold $N$  and a map  $f:M\to BSpin$ such that $f_{|N}$ is constant with value $*$, and which is refined  by a normal $B$-structure. We assume that the Riemannian metric has a product structure near $N$.
Our reason for considering this more general situation where $N$ is not necessarily the boundary of $M$ is that we want to include a case where Lemma \ref{dwedewldew} below  gives a non-trivial result.

\bigskip

 We can assume that $\hat f$ in \eqref{hhdwudiwdwdw} factors over  a map $\hat f_{Spin}:M\to BSpin(k)$ for some $k\in \nat$, and that $\hat f_{Spin|N}$ is constant. We let $P\to M$ be the $Spin(k)$-principal bundle classified by $\hat f_{Spin}$ and form the associated $Spin^{c}(k)$ principal bundle
$\tilde P:=P\times_{Spin(k)}Spin^{c}(k)$. We choose a tangential $Spin^{c}$-structure, i.e.  a $Spin^{c}$-structure $Q\in Spin^{c}(TM)$ and an isomorphism
\begin{equation}\label{fhhefzueuif}Q\otimes \tilde P\cong \underline{Spin^{c}(n+k)}_{M}\ , \end{equation}
 which refines \eqref{udiuqdhqwdwqd}.
 The connection $\tilde \nabla^{TM}$ is  an extension of the Levi-Civita connection to
$Q$.
The bundle $P_{|N}$ is trivialized.
We choose a connection $\nabla^{P}$ on $P$ which restricts to  the trivial connection on $N$. It  further induces a connection $\nabla^{\tilde P}$ on $\tilde P$. 

\bigskip

If $(\theta,V_{\theta})$ is a complex, finite-dimensional representation of $Spin(k)$, then we can define a geometric bundle $\bP(\theta)$ by forming the associated bundle $ P(\theta):=P\times_{Spin(k)}V_{\theta}$ with the induced connection $\nabla^{P(\theta)} $.   If $\iota:V_{\theta}\stackrel{\sim}{\to} \C^{m}$ is an isomorphism of Hilbert spaces, then using the trivialization of $P_{|N}$ we get an isomorphism of geometric bundles $P(\iota):\bP(\theta)_{|N}\to \underline{\C^{m}}_{N}$.  We define
the $\Z/2\Z$-graded bundle $\tilde \bP(\theta)$ such that its even part is $\bP(\theta)$, and  its odd part is $ \underline{\C^{m}}_{M}$. By Lemma \ref{edowedewd} we then get a class
$$\widehat{[\tilde \bP(\theta),P(\iota)]}\in \hat K^{0}(M,N)\ .$$

In the following we show that the class $\widehat{[\tilde \bP(\theta),P(\iota)]}$ depends on $\iota$ in a non-trivial way. Because of this the construction of geometrization along the lines of  \cite[Prop. 5.13]{2011arXiv1103.4217B} has to be modified as will be explained below.  Assume that we have chosen a second isomorphism $\iota^{\prime}$. Then we can write $\iota^{\prime}=\exp(L)\circ  \iota$ for some Lie algebra element $L\in u(\dim(V_{\theta}))$. 
Let $\PD[N]\in H^{1}(M,N;\R)$ denote the dual class of the orientation class $[N]\in H_{n-1}(M;\R)$. 
 \begin{lem}\label{dwedewldew}
We have
$$\widehat{[\tilde \bP(\theta),P(\iota^{\prime})]}-\widehat{[\tilde \bP(\theta),P(\iota)]}=a(\frac{\Tr(L)}{2\pi i}b\ \PD[N])\ .$$
\end{lem}
\proof
We are going to use the homotopy formula.
On
$[0,1]\times M$ we define the  $\Z/2\Z$-graded   bundle $$\hat P:=\pr_{M}^{*}  P(\theta)\oplus \underline{\C^{m}}_{[0,1]\times M}\ .$$ On its restriction to $N$ we consider the isomorphism
$\rho(t):=\exp(tL)\circ \pr^{*}_{N}P(\iota)$, where $t$ is the coordinate of the interval. 
Then $\rho$ interpolates between $P(\iota)$ and $P(\iota^{\prime})$. In order to turn 
$\hat P$ into a geometric bundle $\hat \bP$ we equip its even part  with the connection $\pr_{M}^{*}\nabla^{P(\theta)}$ and the odd part  with the connection 
   $$ \nabla^{triv}-\chi(r)  L dt\ ,$$ where $r:M\to [0,1)$ is the normal coordinate near $N$, $\chi(r)$ is a cut-off function which is equal to $1$ near $r=0$ and vanishes for $r>1/2$.   Then
$\rho: \hat \bP^{+}_{|[0,1]\times N}\to  \hat \bP^{-}_{|[0,1]\times N}$ is an isomorphism of geometric bundles. We have
$$R(\widehat{[ \hat \bP ,\rho]})=\pr^{*}_{M}R \widehat{[\tilde \bP(\theta),P(\iota)]}+
dt\wedge d\chi(r) \frac{\Tr(L)}{2\pi i}b \ .$$
By the homotopy formula for differential $K$-theory \textcolor{black}{(the proof of \cite[Lemma 5.1]{MR2608479}) goes through in the relative case)} we have
$$\widehat{[\tilde \bP(\theta),P(\iota^{\prime})]}-\widehat{[\tilde \bP(\theta),P(\iota)]}=a([d\chi(r)\frac{\Tr(L)}{2\pi i}b,0])\ ,$$
where we use the mapping cone notation $[-,-]$ for forms which was  introduced before Lemma  \eqref{qdlqwd}.
Finally note that $[d\chi(r),0]\in H^{1}(M,N;\R)$ is the Poincar\'e dual class of $[N]\in H_{n-1}(M);\R)$. 
\hB

\begin{kor}\label{cjkdcascsc}
If $N$ is the boundary of $M$, then the class
 $\widehat{[\tilde \bP(\theta),P(\iota)]}$
does not depend on the choice of $\iota$.
\end{kor}
\proof
We have $[N]=0$. \hB

Let  $\tilde R(Spin(k))\subset R(Spin(k))$ be the ideal  of the representation ring of $Spin(k)$ of elements with  vanishing dimension. The associated bundle construction induces a homomorphism
$$\asss:\tilde R(Spin(k))  \to K^{0}(BSpin(k),*)\ .$$ It follows from the completion theorem  \cite{MR0259946} that $\asss$ is injective and has a dense range. 

We choose a basis $(\theta_{i})_{i\in I}$ of the free $\Z$-module $\tilde R(Spin(k))$. We consider the element $\theta_{i}\in \tilde R(BSpin(k))$ as a $\Z/2\Z$-graded representation of $BSpin(k)$. For each 
$i\in I$ we further choose a Hilbert space isomorphism $\iota_{i}$ between the even and odd parts of $\theta_{i}$.
We can define a continuous map
$$\cG_{0}:K^{0}(BSpin(k),*)\to \hat K^{0}(M,N)$$
by the prescription
\begin{equation}\label{cdklncldcd}\cG_{0}(\asss(\theta_{i})):=\widehat{[ \tilde \bP(\theta_{i}), P(\iota_{i})]} \end{equation}
for all $i\in I$.
Then we clearly have the identity
$I\circ \cG_{0}=f^{*}$.
The cohomological character of $\cG_{0}$ is fixed by 
 \begin{equation}\label{ffrfrvffv}c_{\cG_{0}}(\Td^{-1}\cup \ch( \asss(\theta_{i}))):=\Td(\tilde\nabla^{TM})\wedge \ch(\nabla^{\tilde P(\theta_{i})})\in \Omega^{0}_{cl}(M,N)\ .\end{equation}  The map $c_{\cG_{0}}$ only preserves the $b$-degree if the equality of Todd forms $\Td(\tilde\nabla^{TM})=\Td(  \nabla^{\tilde P})^{-1}$ holds. In general this is not the case, and in order to turn $\cG_{0}$ into a geometrization,  we must  add a correction term.
Using the tangential $Spin^{c}$-structure \eqref{fhhefzueuif}   we can define the transgression
\begin{equation}\label{djwkdwd}\delta:=\tilde \Td(\tilde \nabla^{TM}\oplus   \nabla^{\tilde P},\nabla^{triv})\in \Omega P^{-1}(M)/\im(d)\ .\end{equation}
 Then
$$d  (\delta \wedge \Td(    \nabla^{\tilde P})^{-1})=\Td(\tilde \nabla^{TM})-\Td(  \tilde \nabla^{P})^{-1}\ .$$
Let $i:BSpin(k)\to  BSpin$ be the inclusion.
We define the continuous map
$$\cG:K^{0}(BSpin,*)\to \hat K^{0}(M,N)$$ by 
\begin{equation}\label{dkiwoldwdwd} \cG(\phi):=\cG_{0}(i^{*}\phi)-\textcolor{black}{a\left(\frac{\delta\wedge c_{\cG_{0}}(\Td^{-1}\cup \ch(i^{*}\phi))}{\Td(\tilde \nabla^{TM})^{2}} \right)} \ .\end{equation}
\textcolor{black}{We have
\begin{eqnarray*}\lefteqn{\Td(\tilde \nabla^{TM})\wedge R(\cG(\asss(\theta_{i})))}&&\\&=& c_{\cG_{0}}(\Td^{-1}\cup \ch(\asss(\theta_{i})))
-\Td(\tilde \nabla^{TM})\wedge d\left(\frac{\delta\wedge c_{\cG_{0}}(\Td^{-1}\cup \ch(i^{*}\theta_{i}))}{\Td(\tilde \nabla^{TM})^{2}} \right)  \\&=&\Td(\tilde\nabla^{TM})\wedge \ch(\nabla^{\tilde P(\theta_{i})})-\frac{ (\Td(\tilde \nabla^{TM})-\Td(  \tilde \nabla^{P})^{-1})\wedge \Td(\tilde\nabla^{TM})\wedge \ch(\nabla^{\tilde P(\theta_{i})}}{\Td(\tilde \nabla^{TM})})\\
&=&\Td(  \tilde \nabla^{P})^{-1}\wedge \ch(\nabla^{\tilde P(\theta_{i})})\ .
\end{eqnarray*}
Consequently, for the cohomological character of $\cG$ we can take 
the map determined by
$$c_{\cG}(\ch(\asss(\theta))):=\ch(\nabla^{\tilde P(\theta)})\ , \quad \theta\in \tilde R(Spin(k))\ .$$
This homomorphism clearly preserves the $b$-degree.}
Therefore $\cG$ is a geometrization.

In contrast to the absolute case this geometrization not only depends on  $\nabla^{P}$, but also on the choice of the isomorphisms $\iota_{i}$. Nevertheless the construction is sufficiently canonical so that if $(M,N,f,\tilde \nabla^{TM})$ is obtained by taking the boundary $\partial_{B}$ of a zero bordism $(W,F,\tilde \nabla^{TW})$, then
the geometrization extends to $W$.

%
%
%


\section{The intrinsic formula}\label{jcknkwjebckewc}

From now on we consider the notation as in Section \ref{jdkqdqd}. We are going to express the quantity \eqref{jkehehfwiefewff} solely  in terms of data on   a cycle $(M,N,f)$ 
for $x\in \pi_{n}(M(B,A))_{tors}$ such that $\ell x=0$.  Recall that we have fixed  
 a tangential $Spin^{c}$-structure and a $Spin^{c}$-extension  $\tilde \nabla^{TM}$
of the Levi-Civita connection. Let $\cG$ be  a geometrization of $(M,N,f, \tilde \nabla^{TM})$ as   in Definition \ref{iuolqdqdwd}. Let $\phi\in K^{0}(B,A)$ und $(\bV,\rho)$ be as   in Section \ref{jdkqdqd}  such that
$[\bV,\rho]=f^{*}\phi$. Then by \eqref{dwdlwd} there exists a class 
\begin{equation}\label{ejk23jekl32ejl23e}\gamma_{\phi}:=[\alpha_{\phi},\beta_{\phi}]\in \frac{\Omega P^{-1} (M)\oplus \Omega P^{-2}(N)}{\im(\ch)} \ , \quad \alpha_{\phi|N}=d\beta_{\phi} \end{equation}     such that 
\begin{equation}\label{hdhwduwdzuwd}\widehat{[\bV,\rho]}=  \cG(\phi)-a(\gamma_{\phi})\ .\end{equation}
Here we abuse the $[-,-]$-notation for elements in $A^{0}$ and use it in order to write elements in $A^{0}/\im(\ch)$.

\begin{ass}
We assume that there exists a geometrization $\cG_{W,\partial_{A}W}$ which induces the  geometrization  $\cG$ on the $\ell$ copies of $(M,N)$ in the boundary of $\partial_{B}W$.
\end{ass}

\begin{rem}{\rm This assumption is non-trivial. In general not every geometrization on $(M,N)$ can be obtained as such a restriction.
We refer to \cite{2011arXiv1103.4217B} for a detailed discussion. But note that the example of a geometrization constructed in Section \ref{sjkfsfsfsrfsrf} has this property.
} \hB \end{rem}

For a Dirac operator $\Dirac$ we define
  the reduced $\eta$-invariant by
\begin{equation}\label{hgggghj32r897}\xi(\Dirac):=\left[\frac{\eta(\Dirac)+\dim(\ker(\Dirac))}{2}\right]\in \R/\Z\ .\end{equation}

Let   $\Dirac_{ \cD M}\otimes \cD(\bV,\rho)  $ be
the $Spin^{c}$-Dirac operator on the double $\cD M$ twisted by the double of the bundle $(\bV,\rho)$.
\begin{theorem}[Intrinsic formula]\label{widoqdwqd}  The element 
$\eta^{an}(x)\in Q_{n}(M(B,A))$ is represented by the homomorphism
$$K^{0}(B,A)\ni \phi \mapsto [- \int_{M} \Td(\tilde \nabla^{TM})\wedge \alpha_{\phi} -  \int_{ N} \Td(\tilde \nabla^{TN})\wedge \beta_{\phi}]_{\R/\Z}  -\xi(\Dirac_{ \cD M}\otimes\cD(\bV,\rho))\in \R/\Z
$$
\end{theorem}
\proof
First note that the first term on the right-hand side is well-defined since for
$(\alpha_{\phi},\beta_{\phi})\in \im(\ch)$  the sum of the two integrals yields an integer.
We start with the APS index theorem \cite{MR0397797} (compare with \eqref{gdhwdgwdwd}):
$$[\frac{1}{\ell} \ind(\Dirac_{\cD W}\otimes \cD(\bU,\sigma))_{APS}]=[\frac{1}{\ell}\int_{\cD W} \Td(\tilde \nabla^{T\cD W})\wedge \ch(\nabla^{\cD(\bU,\sigma)})]-\xi(\Dirac_{ \cD M}\otimes \cD(\bV,\rho) )\ .$$
 Using the odd $\Z/2\Z$-symmetry of $\cD(\bU,\sigma)$ on $W^{op}\subset \cD W$,  \eqref{hdhwduwdzuwd},  \eqref{jhdjkaw} and Stokes' theorem together with the second equality in \eqref{ejk23jekl32ejl23e}  we calculate
 \begin{eqnarray*}\lefteqn{\int_{\cD W} \Td(\tilde \nabla^{T\cD W})\wedge \ch(\nabla^{\cD(\bU,\sigma)})}&&\\&=&\int_{ W} \Td(\tilde \nabla^{TW})\wedge R(\cG_{W,\partial_{A}W}(\phi)) -\ell \int_{M} \Td(\tilde \nabla^{TM})\wedge \alpha_{\phi} - \ell \int_{ N} \Td(\tilde \nabla^{\partial_{A}W})\wedge \beta_{\phi}\end{eqnarray*}

We now observe that the homomorphism
$$K^{0}(B,A)\ni \phi\mapsto [\frac{1}{\ell}\int_{ W} \Td(\tilde \nabla^{TW})\wedge R(\cG_{W,\partial_{A}W}(\phi))]\in \R/\Z$$
factorizes over the cohomological character $c_{\cG_{W,\partial_{A}W}}$. In view of \eqref{suiwsws} it therefore belongs to the subgroup $\tilde U^{\R}_{n}$.
We conclude that
$\eta^{an}(x)$ is represented by the asserted map. \hB 

\section{Tertiary invariants}\label{ascsc}

In this section we describe the construction of an invariant $\kappa^{top}$ which is a secondary version of the universal $\eta$-invariant and may detect elements in the homotopy of $MA$ which become trivial when mapped to the homotopy of $MB$.  Special cases will be discussed in the subsequent Sections \ref{dsvdsvdsvvvv} and \ref{dioqd}.

We let $n$ be an even integer. In order to simplify matters  we make   the assumption   that $\pi_{n-1}(M(B,A))$ is a torsion group. 
We consider the diagram 
\begin{equation}\label{kkdduqidwd}\xymatrix{\pi_{n-1}(MB)\ar@{.>}[dr]^{\alpha}\ar[r]& \pi_{n-1}(M(B,A))\ar@{-->}[dr]^{can}\ar[d]^{\eta^{top}}\ar[r]&\pi_{n-2}(MA)\ar[r]^{i_{*}} & \pi_{n-2}(MB)\\
&Q^{\R}_{n-1}(M(B,A))\ar[d]^{p}&\ker(i_{*})\ar@{-->}[dl]^{\kappa^{top}}\ar[u]&\\
&Q^{\R}_{n-1}(B,A)&&}\ .\end{equation}
The upper line is a segment of the long exact sequence associated to the fibre sequence \eqref{jkklwdqwd}. \textcolor{black}{This gives the homomorphism denoted by $can$.}
 The map $\alpha$ is the obvious composition.  We define the abelian group
$$Q^{\R}_{n-1}(B,A):= Q^{\R}_{n-1}(M(B,A))/\im(\alpha)\ .$$ Finally, we define the homomorphism $\kappa^{top}$ by the following diagram chase.
Consider an element $y\in \ker(i_{*})$. Then we choose a lift $x\in \pi_{n-1}(M(B,A))$ \textcolor{black}{under $can$}. By assumption it is torsion and therefore in the domain of the universal $\eta$-invariant. By construction,
the image  $p(\eta^{top}(x))\in Q^{\R}_{n-1}(B,A)$ is independent of the choice of the lift.
\begin{ddd}\label{fwfwefwefeweeeee} We define the map 
\begin{equation}\label{lqkdqwpdqwd}\kappa^{top}:\ker(i_{*})\to Q_{n-1}^{\R}(B,A) \end{equation} 
such that $\kappa^{top}(y):=p(\eta^{top}(x))$.
\end{ddd}

\begin{rem}{\rm 
We consider the universal $\eta$-invariant  $\eta^{top}:\pi_{n-1}(MA)_{tors}\to Q_{n-1}^{\R}(MA)$ as a secondary invariant of the $K$-orientation $MA\to BSpin^{c} \stackrel{ABS}{\to} K$. In this sense $\kappa^{top}$ is a tertiary invariant.}\hB 
\end{rem}



\section{Laures' $f$-invariant}\label{dsvdsvdsvvvv}

In this section we discuss an example for the tertiary invariant  defined in \ref{fwfwefwefeweeeee} which has  already  been studied intensively.  We consider the case where 
  $$A=*\ , \quad B= BU\ ,$$ and where $BU\to  BSpin^{c}$ is the canonical map. Then $MA\simeq S$ is the sphere spectrum,  $MU\simeq MB$,  and the corresponding cohomology theories are called framed bordism and  complex bordism. In particular,  the teritary invariant    detects elements in the stable homotopy groups of the sphere.  
   The usual notation for the relative bordism spectrum is \begin{equation}\label{kslxsxisws}\overline{MU}:=M(BU,*)\ .\end{equation}

 The main problem is to define a map out of the group $Q^{\R}_{n-1}(BU,*)$ which is able to detect interesting  elements.  
  
  \bigskip

The construction of the desired evaluation on $Q^{\R}_{n-1}(BU,*)$ employs an elliptic cohomology theory. To this end we fix
an integer   $D\in \nat$ with \textcolor{black}{$
D\ge 4$}  and a $D$'th root of unity $\zeta_{D}$.
Furthermore, we choose a cusp $c$ for the group $\Gamma_{1}(D)$  which is not the cusp at $\infty$.
Then there exists a Landweber exact elliptic cohomology theory  over the ring of  modular forms $$\MF_{*}^{E}:=\textcolor{black}{{}^{c}}\MF_{*}^{\Gamma_{1}(D)}[D^{-1},\zeta_{D}^{-1}]$$ for the group $\Gamma_{1}(D)$
\textcolor{black}{which are holomorphic at all cusps except possibly at $c$}, and  whose $q$-expansions have coefficients in the ring $\Z[D^{-1},\zeta_{D}^{-1}]$ (see e.g. \cite[Thm. 1.2.1]{MR1660325}).
This cohomology theory is represented by a spectrum  $E$ which fits into a sequence of maps
\begin{equation}\label{fmrfrf}MU\to  E\to K[D^{-1},\zeta_{D}][[q]] \ ,\end{equation}  where $MU\to E$ is the complex orientation of the cohomology theory $E$. \textcolor{black}{For a construction of  the map 
 $E\to K[D^{-1},\zeta_{D}][[q]]$  and further details and references we refer to \cite[Sec. 4]{MR2652438}. }  

\bigskip

From now on we assume that $n$   is an
  even integer.   The space of $q$-expansions  $$\MF^{E}_{n/2}[[q]]\subseteq \Z[D^{-1},\zeta_{D}][[q]] \cong  \pi_{n}(K[D^{-1},\zeta_{D}][[q]])$$  is  the image of $\pi_{n}(E)$ under this evaluation. 

\bigskip

We extend the composition \eqref{fmrfrf}   
 to a composition of maps of vertical fibre sequences
\begin{equation}\label{jhbewhjkwed}\xymatrix{S\ar[d]\ar@{=}[r]&S\ar[d]\ar[r]^{\varepsilon}&\ar[d]K[D^{-1},\zeta_{D}]  \\MU\ar[d]\ar[r]
&E\ar[r]\ar[d]& K[D^{-1},\zeta_{D}][[q]]\ar[d] \\\overline{MU}\ar[r]\ar@/_1cm/[rrr]_{\phi}&
\overline{E}\ar[r]&  K[D^{-1},\zeta_{D}][[q]]/K[D^{-1},\zeta_{D}]\ar[r]^{\simeq}&\prod_{i=1}^{\infty}  q^{i}K[D^{-1},\zeta_{D}]}\ .\end{equation}
The map $\phi$ is defined as the natural composition.
We interpret $\phi$ as a sequence of classes $\phi_{i}\in K[D^{-1},\zeta_{D}]^{0}(\overline{MU})$ defined  for all positive $i\in \nat$.
 For even $n$ the evaluation at $\phi$ induces a map $$\widetilde{\ev}_{\phi}:\Hom^{cont}(K^{0}(\overline{MU}),\pi_{n}(K\R/\Z))\to \frac{\C[[q]]}{\Z[D^{-1},\zeta_{D}][[q]]+q^{0}\C+\MF_{\C,n/2}^{E}[[q]]}\ , $$ $$\widetilde{\ev}_{\phi}(h):=\left[\sum_{i=1}^{\infty} \ev_{\phi_{i}}(h) q^{i}\right]\ . $$
In this formula we interpret $\ev_{\phi_{i}}(h) \in \pi_{n}(K\R/\Z)\cong \R/\Z$ and use the well-defined homomorphisms
$$\R/\Z\to \frac{\C[[q]]}{\Z[D^{-1},\zeta_{D}][[q]]+q^{0}\C+\MF_{\C,n/2}^{E}[[q]]}\ , \quad [x]\mapsto [xq^{i}]\ , \quad x\in \R$$ for $i\in \nat$, $i\ge 1$. 
\textcolor{black}{Finally, we have set $$\MF_{\C,n/2}^{E}[[q]]:= \MF_{n/2}^{E}[[q]]\otimes_{\Z[D^{-1},\zeta_{D}]} \C\subset \C[[q]]\ .$$}
 
 We now observe that for $y\in U_{n-1}^{\R}$ we have $\widetilde{\ev}_{\phi}(y)\in \MF_{\C,n/2}^{E}[[q]]+q^{0}\C$. Therefore $\widetilde \ev_{\phi}$ descends to  a homomorphism
$$ \ev_{\phi}:Q^{\R}_{n-1}(\overline{MU})\to  \frac{\C[[q]]}{\Z[D^{-1},\zeta_{D}][[q]]+ q^{0}\C  + \MF_{\C,n/2}^{E}[[q]]}\ . $$

 Since $n$ is even we have $\pi_{n-1}(MU)\cong 0$ and there are no interesting elements which can be detected by the evaluation of the universal $\eta$-invariant for $MU$ using this evaluation. On the other hand, this fact implies that the evaluation $\ev_{\phi}$
 actually further  descends to a  homomorphism
$$\overline{ \ev}_{\phi}:Q^{\R}_{n-1}(BU,*)\to \frac{\C[[q]]}{\Z[D^{-1},\zeta_{D}][[q]]+ q^{0}\C  + \MF_{\C,n/2}^{E}[[q]]}\ .$$

We now assume in addition  that $n$ satisfies $n\not=2$. Then we have the equality
  $$\ker\left(i_{*}:\pi_{n-2}(S)\to \pi_{n-2}(MU)\right)=\pi_{n-2}(S)\ .$$ 
  \begin{ddd}  For $n\not=2$ the $f$-invariant is defined to be the homomorphism  
$$f:=\overline \ev_{\phi}\circ \kappa_{MU}^{top}:\pi_{n-2}(S)\to \frac{\C[[q]]}{\Z[D^{-1} ,\zeta_{D}][[q]]+ q^{0}\C  + \MF_{\C,n/2}^{E}[[q]]}\ . $$
\end{ddd} 
Here $\kappa^{top}_{MU}$ is the specialization of $\kappa^{top}$ in \eqref{kkdduqidwd} to the present case. The name $f$-invariant is justified by the fact 
 verified in \cite{MR2652438}  that $f$ is indeed  Laures' $f$ invariant introduced in  \cite{MR1660325}. For explicit calculations we refer to  \cite{MR1660325}, 
\cite{MR1781277} and \cite{2008arXiv0808.0428V}.
In particular, the $f$-invariant is non-trivial.

\bigskip

\textcolor{black}{Some properties of the $f$-invariant can easily be deduced from properties of the universal $\eta$-invariant shown in \cite[Prop. 2.7]{2011arXiv1103.4217B}. We will explain this in the following. Using much deeper results of  \cite{MR1660325} we will also see that $\kappa^{top}_{MU}$ is non-trivial.}

The relative bordism spectrum $\overline{MU}$ is the main constituent of the Adams tower
 \begin{equation}\label{wehjklwed}S\leftarrow \Sigma^{-1}\overline{MU}\leftarrow \Sigma^{-1}\overline{MU}\wedge  \Sigma^{-1} \overline{MU}\leftarrow\Sigma^{-1}\overline{MU}\wedge  \Sigma^{-1} \overline{MU}\wedge  \Sigma^{-1} \overline{MU}\leftarrow\dots\ .\end{equation}
 The Adams tower induces a decreasing filtration $F^{*}_{MU}\pi_{m}(W)$  of the  $m$th homotopy group of a spectrum $W$ for every $m\in \Z$. In order to define this filtration we consider the smash product of the tower \eqref{wehjklwed}  with $W$.
 For $k\in \nat$ we define $F^{k}_{MU}\pi_{m}(W)\subseteq \pi_{m}(W)$ to be the subgroup of elements   which lift to $W\wedge \Sigma^{-k}\overline{MU}^{\wedge k}$.

For all integers $k\ge 0$ and even $n$ we have the obvious isomorphisms
$$F^{k}_{MU}\pi_{n-2}(S)\cong F_{MU}^{k-1}\pi_{n-1}(\overline{MU})\ .$$
In particular for $n\not=2$ we have
$$\pi_{n-2}(S)\cong F^{2}_{MU}\pi_{n-2}(S)\cong F^{1}_{MU}\pi_{n-1}(\overline{MU})\ .$$
The complex orientation $MU\to K$ of $K$-theory induces a map
$$F^{3}_{MU} \pi_{n-2}(S)\cong F^{2}_{MU}\pi_{n-1}(\overline{MU})\to F^{2}_{K}\pi_{n-1}(\overline{MU})\ ,$$ where the filtration $F^{*}_{K}$ is defined similarly replacing $MU$ by $K$ in the construction of the Adams tower.
By \cite[Prop. 2.7]{2011arXiv1103.4217B} the image of this map is annihilated by $\eta^{top}$.
Since $K_{*}(\overline{MU})$ is torsion-free, \textcolor{black}{by \cite[Prop. 2.7, 4]{2011arXiv1103.4217B} the homomorphism $$\bar \eta^{top}:\Gr^{2}_{K}\pi_{n-2}(S)\cong \Gr^{1}_{K}\pi_{n-1}(\overline{MU})\to Q^{\R}_{n-1}(\overline{MU})\cong Q^{\R}_{n-1}(BU,*)$$ induced by $\eta^{top}$ is injective. 
Furthermore, for $\kappa_{MU}^{top}$ we get a factorization 
\begin{equation}\label{jbdjkewbdkjwedewde}\bar \kappa^{top}_{MU}:\Gr^{2}_{MU}\pi_{n-2}(S)\stackrel{!}{\to} \Gr^{2}_{K}\pi_{n-2}(S)\stackrel{\bar\eta^{top}}{\to} Q^{\R}_{n-1}(BU,*)\ .\end{equation}  Consequently, our theory implies that the $f$-invariant induces  a homomorphism}
\begin{equation}\label{kdqdlqdqwd}\bar f: \Gr^{2}_{MU}\pi_{n-2}(S)\textcolor{black}{\stackrel{\bar \kappa^{top}_{MU}}{\longrightarrow}} Q^{\R}_{n-1}(BU,*) \textcolor{black}{\stackrel{\bar \ev_{\phi}}{\longrightarrow}}  \frac{\C[[q]]}{\Z[D^{-1},\zeta_{D}][[q]]+ q^{0}\C  + \MF_{\C,n/2}^{E}[[q]]}\ .\end{equation}

\textcolor{black}{
While  \cite[Prop. 2.7]{2011arXiv1103.4217B}  implies that the $f$-invariant factorizes over $\Gr_{MU}^{2}\pi_{n-2}(S)$ we need the much deeper result of 
  \cite{MR1660325} in order to see that it is non-trivial. In this paper it has been shown that 
  the homomorphism \eqref{kdqdlqdqwd} is injective.
  Consequently,
$\overline{\ev}_{\phi}$ detects the image of $\bar \kappa^{top}_{MU}$.
Moreover,  we conclude that $\bar \kappa^{top}_{MU}$ is non-trivial since it detects the group
$\Gr_{MU}^{2}\pi_{n-2}(S)$ which is non-trivial for many $n$. Finally, we can conclude that the marked homomorphism in 
\eqref{jbdjkewbdkjwedewde} is injective.}

\begin{rem}{\rm 
In \cite{MR2652438} we gave an intrinsic formula for the $f$-invariant in terms of a sequence of $\eta$-invariants of Dirac operators twisted by bundles derived from the tangent bundle. In order to interpret this as an example of  the intrinsic formula \ref{widoqdwqd} it would be necessary to translate
this construction to a construction with a complementary bundle. This would require to extend the theory of geometrizations to $K$- and differential $K$-theory with coefficients in $\Z[D^{-1},\zeta_{D}]$. We think that this is possible but that the details are not very enlightning. We will demonstrate the intrinsic formula in a second example in Section \ref{dioqd}.} \hB 
  \end{rem}

\section{A $Spin$-version of the $f$-invariant}\label{dioqd}

In this section we discuss a $Spin$-version $\kappa^{top}_{MSpin}$ of the tertiary invariant
defined in \ref{fwfwefwefeweeeee}. We will see in Proposition  \ref{dedlwd} that it is non-trivial.
Similarly as in Section \ref{dsvdsvdsvvvv} we define an evaluation leading to a $Spin$-version $f^{Spin}$ of the $f$-invariant. We make the intrinsic formula for $f^{Spin}$ explicit.

\bigskip
 
We consider the case where  
 $$A=*\ , \quad  B=BSpin\ ,$$ and we choose the canonical map $BSpin\to BSpin^{c}$. In analogy to \eqref{kslxsxisws} we set  $$\overline{MSpin}:=M(BSpin,*)\ .$$ In the following we construct a sequence of $K$-theory classes $\phi_{i}\in K^{0}(BSpin,*)$, $i\ge 1$, which will be used to define an evaluation on $Q^{\R}_{n-1}(BSpin,*)$.

\textcolor{black}{In the following
a complex vector bundle will considered it as $\Z/2\Z$-graded such that
its even part is the given bundle and its odd part is zero. For a $\Z/2\Z$-graded complex vector 
$E$ we write $E^{\pm}$ for its even and odd components, and $E^{op}$ for the bundle with the opposite grading, i.e. with the even and odd parts exchanged. For a complex vector space $W$  and a manifold $M$ we write $\underline{W}_{M}$ for the trivial  complex vector bundle with fibre $W$ on $M$.}
 
For a  real vector bundle $V$ on a manifold or space $M$  we consider the $\Z/2\Z$-graded complex vector bundle $$\tilde V:=V\otimes \C\oplus  (\underline{\C^{\dim(M)}}_{M})^{op}$$ and define the formal power series of $\Z/2\Z$-graded complex vector bundles
\begin{equation}\label{gdhjdgwqgdjdqdqwd}\Theta(V ):=\bigotimes_{u=1}^{\infty}S_{q^{u}}(\tilde V) \otimes\bigotimes_{v=1}^{\infty}\otimes \Lambda_{-q^{v-\frac{1}{2}}}(\tilde V)\ .\end{equation}
To be precise, in this definition we  expand the tensor products and sort the terms according to the powers of $q^{1/2}$. In particular,  we  interpret a summand $-W$ as $W^{op}$. In this way we get a series of complex vector bundles
$$\Theta(V )=\sum_{i=0}^{\infty} q^{i/2} \Theta_{i}(V)\ ,$$
where the $\Z/2\Z$-graded complex vector bundles $\Theta_{i}(V)$ are constructed from $V$ in a functorial way using operations of the tensor calculus. 
The power series is multiplicative in the sense that $$\Theta(V\oplus V^{\prime})\cong \Theta(V )\otimes \Theta(V^{\prime} )\ .$$
Furthermore, since
$$[\Theta(V )]=1+O(q^{1/2})\in K[[q^{1/2}]]^{0}(M)\ ,$$ the $K$-theory
class  $[\Theta(V )]\in K[[q^{1/2}]]^{0}(M)$ is a multiplicative unit.
Therefore we can extend the association $V\mapsto [\Theta(V)]$ to real $K$-theory classes, i.e. we can define a $K$-theory class $\Theta(\xi )\in K[[q^{1/2}]]^{0}(M)$ for a $K$-theory class  $\xi\in KO^{0}(M)$.

\bigskip

Since the construction  $V\mapsto \Theta_{i}(V)$ is functorial,
a trivialization $\rho:V\stackrel{\sim}{\to} \underline{\R^{k}}$ induces a sequence of isomorphisms $$\Theta_{i}(\rho) :
\Theta_{i}(V)^{+}\to \Theta_{i}(V)^{-}\ , \quad i\ge 1$$
\textcolor{black}{(recall that $\Theta_{i}(V)^{\pm}$ denote the even and odd parts of the $\Z/2\Z$-graded bundle $\Theta_{i}(V)$).}

 Let  $\xi\in KO^{0}(BSpin)$ be the class of the normalized  universal bundle.
Then we get a class
$$\Theta(\xi ) \in K[[q^{1/2}]]^{0}(BSpin)\ .$$

We have a preferred trivialization $\iota$ of $\xi_{|*}$ so that we get
a sequence of relative classes
$$[\Theta_{i}(\xi),\Theta_{i}(\iota)]\in K^{0}(BSpin,*)\ , \quad i\ge 1\ .$$
The series $$1+\sum_{i=1}^{\infty}q^{i/2}[\Theta_{i}(\xi),\Theta_{i}(\iota)]\in 1+q^{1/2}K[[q^{1/2}]](BSpin,*)$$ is invertible. 
We define a sequence of classes $\phi_{i}\in K[[q^{1/2}]](BSpin,*)$ indexed by positive integers $i$   uniquely such that
\begin{equation}\label{dkqolwdwdqqqq}1+\sum_{i=1}^{\infty}q^{i/2} \phi_{i}=\left(1+\sum_{i=1}^{\infty}q^{i/2}[\Theta_{i}(\xi),\Theta_{i}(\iota)]\right)^{-1}\ .\end{equation}


\bigskip

For $k\in \Z$ let 
 $$\MF^{\Gamma^{0}(2),\R}_{k} [[q^{1/2}]]\subseteq \R[[q^{1/2}]]$$  denote the groups of $q^{1/2}$-expansions of holomorphic   modular forms of weight $k$ for the congruence group $\Gamma^{0}(2)\subseteq SL(2,\Z)$ with real Fourier coefficients at the cusp at $\infty$.

 We now assume that $n\in \Z$ is even. The evaluation at the classes $\phi_{i}$  induces a map $$\widetilde{\ev}_{\phi}:\Hom^{cont}(K^{0}(\overline{MSpin}),\pi_{n}(K\R/\Z))\to \frac{\R[[q^{1/2}]]}{\Z [[q^{1/2}]]+q^{0}\R+\MF^{\Gamma^{0}(2),\R}_{n/2} [[q^{1/2}]]}\ ,$$ $$ \widetilde{\ev}_{\phi}(h):=\left[\sum_{i=1}^{\infty} \ev_{\phi_{i}}(h) q^{i/2}\right]\ . $$
 
 
If $y\in U_{n-1}^{\R}$, then by \cite[Prop. 2.6]{MR2153079} we know that  $\widetilde{\ev}_{\phi}(y)\in \MF^{\Gamma^{0}(2),\R}_{n/2} [[q^{1/2}]]+q^{0}\R$. Therefore $\widetilde \ev_{\phi}$ descends to  a homomorphism
$$\ev_{\phi} :Q^{\R}_{n-1}(\overline{MSpin})\to  \frac{\R[[q^{1/2}]]}{\Z [[q^{1/2}]]+ q^{0}\R  + \MF^{\Gamma^{0}(2),\R}_{n/2} [[q^{1/2}]]} \ .$$
 \begin{lem}\label{fwioelfewf}
If $n\equiv 0 (4)$, then the universal $\eta$-invariant  $$\eta^{top}:\pi_{n-1}(MSpin)\to Q^{\R}_{n-1}(MSpin)$$ is trivial.
\end{lem}
\proof 
One checks using  the calculation of $MSpin$ by \cite{MR0219077} and the results of   \cite{MR0231369} that the Bousfield $K$-localization map  $\pi_{n-1}(MSpin)\to \pi_{n-1}(MSpin_{K})$ is trivial.
The assertion now follows from the fact that $\eta^{top}$ factors through the $K$-localization \cite[Lemma 2.8]{2011arXiv1103.4217B}. \hB 

From now on we assume that $n\equiv 0(4)$. Lemma  \ref{fwioelfewf} implies that
  $\alpha:\pi_{n-1}(MSpin)\to Q_{n-1}^{\R}(\overline{MSpin})$ vanishes, where $\alpha$ is defined in  \eqref{kkdduqidwd}.
Therefore the evaluation $\ev_{\phi}$
 actually factorizes over
$$\overline{ \ev}_{\phi}:Q^{\R}_{n-1}(BSpin,*)\to \frac{\R[[q^{1/2}]]}{\Z [[q^{1/2}]]+ q^{0}\R  + \MF^{\Gamma^{0}(2),\R}_{n/2} [[q^{1/2}]]}\ .$$

\textcolor{black}{
It is known that for $k\in \nat$ the composition
$$\pi_{8k+2}(S)\to \pi_{8k+2}(MSpin)\to \pi_{8k+2}(KO)\cong \Z/2\Z$$ is surjective.
Consequently the domain of definition of $\kappa^{top}_{MSpin}$ in Definition  \ref{fwfwefwefeweeeee}  is the in general proper subgroup $$\ker(i_{*}):=\ker\left(i_{*}: \pi_{n-2}(S)\to \pi_{n-2}(MSpin)\right)\subseteq \pi_{n-2}(S)\ .$$}
 
\begin{ddd} For every integer $n$ with $n\equiv 0(4)$ we define 
the $Spin$-version of the $f$-invariant by 
$$f^{Spin}:=\overline{ \ev}_{\phi}\circ \kappa^{top}_{MSpin} :\textcolor{black}{\ker(i_{*})}\to \frac{\R[[q^{1/2}]]}{\Z [[q^{1/2}]]+ q^{0}\R  + \MF^{\Gamma^{0}(2),\R}_{n/2} [[q^{1/2}]]}\ .$$
\end{ddd}

At the moment we have the following information about $f^{Spin}$.

\begin{prop} \label{dedlwd}\begin{enumerate}\item The $Spin$ version $f^{Spin}$ of the $f$-invariant vanishes $2$-locally.
\item The tertiary invariant $\kappa^{top}_{MSpin}$ is non-trivial and detects, at least, 
the odd torsion in $\Gr^{2}_{MU}\pi_{n-2}(S)$.
\end{enumerate}
\end{prop}
\proof

\bigskip
\textbf{The two-local case}
\bigskip

In the following all spectra a two-localized. We have the ABP-splitting \cite{MR0219077}
$$MSpin\simeq  F\vee G\ ,$$
where $F$ is a wedge of truncated copies of $ko$ and $G$ is a wedge of
shifted copies of $H\Z/2\Z$. We have a factorization of the unit $S\to  MSpin$ as
$$S\to F\to F\vee G\ .$$ For a spectrum $X$  Bousfield localization at $K$ leads to a fibre sequence
$$X^{K}\to X\to X_{K}\to \Sigma X^{K}\ .$$
By  \cite{MR0231369} we have $G_{K}\simeq 0$ and therefore $G^{K}\stackrel{\sim}{\to} G$.
We apply $K$-localization to the sequence
$$\Sigma^{-1}\overline{MSpin}\to S\to MSpin\ .$$

We get the following web of horizontal and vertical fibre sequences
$$\xymatrix{\Sigma^{-1} \overline{MSpin}^{K}\ar[r]\ar[d]&S^{K}\ar[d]\ar[r]&F^{K}\vee G\ar[d]\\
\Sigma^{-1} \overline{MSpin}\ar[r]\ar[d]&S\ar[r]\ar[d]&F\vee G\ar[d]\\
\Sigma^{-1} \overline{MSpin}_{K}\ar[r]& S_{K}\ar[r]&F_{K} }\ .$$
The right upper map  factorizes as $S^{K}\to F^{K}\to F^{K}\vee G$.
Recall that we assume that $n\equiv 0(4)$. We have an exact sequence
$$\pi_{n-1}(F)\to \pi_{n-1}(F_{K})\to \pi_{n-2}(F^{K})\to \pi_{n-2}(F)\stackrel{(i)}{\to} \pi_{n-2}(F_{K})\ .$$
We now use that the canonical map $ko[a..\infty] \to KO$  becomes the $K$-localization map after $2$-completion. In view of the structure of $F$
this implies that the  map $(i)$ is injective, and that $\pi_{n-1}(F_{K})= 0$. We conclude that  $\pi_{n-2}(F^{K})=0$.
 
 \bigskip

 \begin{lem}\label{dnekjdwedw}
 The natural homomorphism  $\ker(i_{*})_{(2)} \to \pi_{n-2}(S_{K})$ vanishes. 
 \end{lem}
 \proof
 Let $KO_{2}$ be $KO$ completed at $2$ and  $g\in \Z^{\times}_{2}$ be a topological generator of
 $\Z^{\times}/\{\pm 1\}$.
 We consider the following commuting diagram:
  \begin{equation}\label{hhjefgjwefewf}\xymatrix{&S\ar[r]^{i}\ar[d]&MSpin\ar[d]^{ABS}& \\
 &S_{K}\ar@{..>}[r]\ar@{-->}[d]&KO\ar[d]^{2-compl}&\\\Sigma^{-1} KO_{2}\ar[r]&S_{K(1)}\ar[r]&KO_{2}\ar[r]^{(\psi_{g}-1)}&KO_{2}}\ .\end{equation}
 The dotted arrow is obtained by the universal property of $K$-localization using that $KO$ is $K$-local. 
 The lower  part is a fibre sequence which defines the $K(1)$-local sphere.
 The dashed arrow is obtained as a lift of the obvious composition $S_{K}\to KO\to KO_{2}$.

 Since $\pi_{4*-1}(KO_{2})=0$ we get  from the associated long exact sequence in homotopy   that $$\pi_{4*-2}(S_{K(1)})\cong \ker\left((\psi_{g}-1):\pi_{4*-2}(KO_{2})\to \pi_{4*-2}(KO_{2})\right)\ .$$ 

Using the arithmetic square
$$\xymatrix{S_{K}\ar[d]\ar[r]&S_{K(1)}\ar[d]\\
S_{H\Q}\ar[r]&(S_{K(1)})_{H\Q}}$$
and the fact that the homotopy of the spectra in the  bottom line vanishes in positive degrees, we see that
the top line induces an isomorphism
$\pi_{n-2}(S_{K})\stackrel{\cong}{\to} \pi_{n-2}(S_{K(1)})$ for $n\ge 2$.
We now use the commutativity of \eqref{hhjefgjwefewf} in order to conclude that for a class $y\in \pi_{n-2}(S)$ the condition $y\in \ker(i_{*})_{(2)}$ implies that its image in $\pi_{n-2}(S_{K})$ vanishes. \hB

 \color{black}

  Let $y\in \ker(i_{*})_{(2)}\subseteq \pi_{n-2}(S)$.
 Then by Lemma \ref{dnekjdwedw} 
 \color{black}  the image of $y$ in  $\pi_{n-2}(S_{K})$ vanishes. We choose a lift   $z\in \pi_{n-2}(S^{K})$. Its image in $\pi_{n-2}(F^{K}\vee G)$ vanishes so that
 we can find a further lift $\tilde z\in \pi_{n-1}(\overline{MSpin}^{K})$. The image $x\in \pi_{n-1}(\overline{MSpin})$ of $\tilde z$ can serve as lift of $y$ in the construction of $\kappa^{top}_{MSpin}$, see Section \ref{ascsc}.
 But since the universal $\eta$-invariant factorizes over the $K$-localization we see that
 $\kappa_{MSpin}^{top}(y)=0$ since the image of $x$ in $\pi_{n-1}( \overline{MSpin}_{K})$ vanishes.

 \bigskip
\textbf{The odd torsion}
\bigskip

We now localize all spectra at an odd prime $p$.
Then we have $MSpin\simeq MSp$. Furthermore, by a result of Baker-Morava \cite{bm}
we know that $MSp$ is a summand of $MU$. It follows that the sequence
$$\Sigma^{-1}\overline{MSpin}\to S \to MSpin$$ is a summand of
the sequence
$$\Sigma^{-1}\overline{MU}\to S\to MU\ .$$ 
\textcolor{black}{This implies that  
$\ker(i_{*})_{(p)}=\pi_{n-2}(S)$ so that, after localization at $p$, the homomorphism $\kappa^{top}_{MSpin}$ is defined on all of $\pi_{n-2}(S)$.
}
 We further get the commutative diagram  
$$\xymatrix{ \pi_{n-2}(S) \ar@/^1cm/[rr]^{\kappa^{top}_{MSpin}}\ar[r]\ar@{=}[d]&\pi_{n-1}(\overline{ MSpin})\ar[r]^{\eta^{top}}\ar@^{(->}[d]& Q^{\R}_{n-1}(\overline{MSpin})_{(p)}\ar[d] \\\pi_{n-2}(S)\ar[r]\ar@/_1cm/[rr]_{\kappa^{top}_{MU}}&\pi_{n-1}(\overline{ MU})\ar[r]^{\eta^{top}} &Q^{\R}_{n-1}(\overline{MU})_{(p)}
}\ ,$$ where the first map is the construction $y\mapsto x$.
Since we know that Laures' $f$-invariant detects $\Gr^{2}_{MU}\pi_{n-2}(S)$ we conclude that 
$$\bar \kappa^{top}_{MSpin}:\Gr^{2}_{MU}\pi_{n-2}(S)\to \frac{Q^{\R}_{n-1}(\overline{MSpin})_{(p)}}{\kappa^{top}_{MSpin}(F_{MU}^{3}\pi_{n-2}(S))}$$
is injective, too.
\hB

We now derive an intrinsic formula for $f^{Spin}$. Let $y\in \pi_{n-2}(S)$ and
consider a cycle   $(M,N,f,\tilde \nabla^{TM})$   for a lift $  x\in \pi_{n-1}(\overline{MSpin})$ of $y$. 
We use the notation introduced  in Section \ref{sjkfsfsfsrfsrf} and assume that 
 $\cG$ is a geometrization as constructed there.

Forming the analog of  \eqref{dkqolwdwdqqqq} on the level of bundles we define a sequence of $\Z/2\Z$-graded vector bundles $\Phi_{i}:=\Phi_{i}(\hat f_{Spin}^{*}\xi^{k})$ and isomorphisms $\iota_{i}  :\Phi_{i|N}^{+}\to \Phi_{i|N}^{-}$ such that  
$$f^{*}\phi_{i}=[\Phi_{i} ,\iota_{i}]\ , \quad i\ge 1\ .$$
Note again, that the bundles $\Phi_{i} $ are obtained from $\hat f_{Spin}^{*}\xi^{k}$ in a functorial way using the operations  of tensor calculus. In particular, the
  connection   $\nabla^{P}$ induces a connection on $\hat f_{Spin}^{*}\xi^{k}$  which is
compatible with the trivialization at $N$, and therefore connections $ \nabla^{\Phi_{i} } $ on $\Phi_{i} $ which are compatible with $\iota_{i} $ for all $i\ge 1$.
We therefore get geometric bundles
$(\mathbf{\Phi}_{i} ,\iota_{i})$,
and by Lemma \ref{edowedewd}  differential $K$-theory classes
$$\hat \phi_{i}:=\widehat{[\mathbf{\Phi}_{i} ,\iota_{i}]}\in \hat K^{0}(M,N)\ .$$
Since $N$ is the boundary of $M$, by Corollary \ref{cjkdcascsc} the classes $\hat \phi_{i}$ do not depend on   $\iota_{i}$ and we thus have in $\hat K^{0}(M,N)$:
$$\hat \phi_{i}=\cG_{0}(\phi_{i})\ .$$

We define the correction forms
$$\gamma_{i}=[\alpha_{i},\beta_{i}]\in \frac{\Omega P^{-1}(M)\oplus \Omega P^{-2}(N)}{\im(\ch)}\ , \quad \textcolor{black}{\alpha_{i|N}=d\beta_{i}}$$  such that
$$ \hat \phi_{i}=\cG(\phi_{i})-a(\gamma_{i}) .$$
From \eqref{dkiwoldwdwd} and \eqref{ffrfrvffv}  we see that we can take
$$\gamma_{i}=[\frac{\delta\wedge c_{\cG_{0}}(\Td^{-1}\cup \ch(\phi_{i}))}{\Td(\tilde \nabla^{TM})^{2}}  ,0 ]=[\frac{\delta\wedge \ch(\nabla^{\Phi_{i} })}{\Td(\tilde \nabla^{TM})},0]\ ,$$
where $\delta$ is defined in \eqref{djwkdwd}.

By specializing Theorem \ref{widoqdwqd} we get the desired intrinsic formula for $f^{Spin}(y)$.
\begin{prop} 
The class $f^{Spin}(y)= \ev_{\phi}(\eta^{top}(x))$ is represented by the series
$$\sum_{i=1}^{\infty} q^{i/2}\left([-\int_{M}\delta\wedge  \ch(\nabla^{\Phi_{i}}) ]-\xi(\Dirac_{\cD M}\otimes \cD(\mathbf{\Phi}_{i} , \iota_{i}))\right)\ . $$
\end{prop}

 \begin{rem}{\rm
The analogy with Laures' $f$-invariant lets us strongly believe, that the $Spin$-version $f^{Spin}$ is non-trivial, too. Since we know that $\kappa_{MSpin}^{top}$ is non-trivial, the remaining  question is whether the evaluation $\overline \ev_{\phi}$ is strong enough to detect some non-trivial elements in the image of $\kappa^{top}_{MSpin}$. The first  case to  check would be the evaluation of $\kappa^{top}_{MSpin}(\beta)$ for  a non-trivial three-torsion element $\beta\in \pi_{10}(S)\cong \Z/6\Z$.  
At the moment the calculation of examples seems to be a non-trivial matter.
\textcolor{black}{In comparison with the case of Laures $f$-invariant a missing ingredient is a homotopy theoretic interpretation of the power series $\Theta$ in terms of a suitable Landweber exact elliptic cohomology theory.}
}\hB
\end{rem}

\begin{rem}{\rm 
This intrinsic formula again works with bundles associated to a geometric normal $Spin^{c}$-structure.
For esthetic reasons it would be interesting to have a formula which uses geometric bundles
associated to the tangent bundle.} \hB
\end{rem}

\bigskip

We finally show how Lemma \ref{fwioelfewf} leads to an improvement of \cite[Thm 1.1]{2013arXiv1312.7494H}. 
\textcolor{black}{In a certain sense this is a side result since  it does not use the extension of the universal eta invariant to the relative case.}
Let $m\in \nat$ be a positive integer.  We consider a closed $4m-1$-dimensional Riemannian spin manifold $M$ with metric $g^{TM}$.  
We equip the tangent bundle $TM$ with the Levi-Civita connection $\nabla^{TM}$ and thus define the  geometric bundle $\mathbf{TM}:=(TM,g^{TM},\nabla^{TM})$. We consider the formal power series (in $q^{1/2}$) of $\Z/2\Z$-graded geometric bundles $ \Theta(\btmr) $, where $\Theta$ is defined in \eqref{gdhjdgwqgdjdqdqwd}.  For $k\in \Z$ \textcolor{black}{we consider the $\R$-vector space  $\MF^{\Gamma^{0}(2),\R,*}_{k}$ of  $\Gamma^{0}(2)$-invariant meromorphic sections  of the $k$-th power of the canonical bundle on the upper half plane which are holomorphic at the cusp at $\infty$ and have real Fourier coefficients there. We furthermore let 
 $$ \MF^{\Gamma^{0}(2),\R,*}_{k} [[q^{1/2}]]\subseteq \R[[q^{1/2}]]$$ denote the group  of $q^{1/2}$-expansions of  the elements of $\MF^{\Gamma^{0}(2),\R,*}_{k}$.}
 Recall the definition \eqref{hgggghj32r897} of the reduced $\eta$-invariant of a Dirac operator.
 We set
$$\xi^{\R}(\Dirac):=\frac{\eta(\Dirac)+\dim(\ker(\Dirac))}{2}\in \R$$
so that $\xi(\Dirac)=[\xi^{\R}(\Dirac)]$.
   The following Theorem was shown in \cite[Thm 1.1]{2013arXiv1312.7494H}.
\begin{theorem}[Han-Zhang] \label{thm1001}
The formal power series   of reduced $\eta$-invariants
$$\xi^{\R}(\Dirac\otimes \Theta(\btmr) )\in \R[[q^{1/2}]]$$
  belongs to the subgroup
$$\Z[[q^{1/2}]]+\MF^{\Gamma^{0}(2),\R,*}_{2m} [[q^{1/2}]]\subseteq \R[[q^{1/2}]]\ .$$  \end{theorem}

Note that $\MF^{\Gamma^{0}(2),\R}_{2m} [[q^{1/2}]]\subseteq \MF^{ \Gamma^{0}(2),\R,*}_{2m} [[q^{1/2}]]$ is the subspace of $q^{1/2}$-expansions of holomorphic modular forms.
Lemma \ref{fwioelfewf} 
 implies:  \begin{theorem}\label{thm1000} The formal power series of reduced $\eta$-invariants
$$\xi^{\R}(\Dirac\otimes \Theta(\btmr))\in \R[[q^{1/2}]]$$
  belongs to the subgroup
$$\Z[[q^{1/2}]]+ \MF^{ \Gamma^{0}(2),\R }_{2m} [[q^{1/2}]] \subseteq \R[[q^{1/2}]] \ .$$  \end{theorem}
 
 \begin{rem}{\rm Theorem \ref{thm1000} improves Theorem \ref{thm1001} since it replaces the condition "meromorphic" by the stronger condition "holomorphic".
 The proof of Theorem \ref{thm1001} given by Han-Zhang depends on a Theorem of Hopkins  which is a consequence of Snaith's theorem and stated in the thesis of Klonoff, 2008. It is very different from the proof of Theorem \ref{thm1000}. 
 
 }\hB \end{rem} 
 
\proof

%
%
%

The given data induces the cycle $(M,  f,\tilde \nabla^{TM})$ for a class $[M,f]\in \pi_{4m-1}(MSpin)$.   By \cite[Prop. 5.13]{2011arXiv1103.4217B} the Levi-Civita connection on $TM$ induces a good geometrization 
$$\cG^{LC}:K^{0}(BSpin)\to \hat K^{0}(M)$$ of  the cycle $(M,f,\tilde\nabla^{TM})$.   It follows from the construction of this geometrization   that 
\begin{equation}\label{eq2c}\cG^{LC}(\Theta(-\xi))=[\Theta(\btmr)  ] \end{equation} 
in $\hat{K}^{0}(M)[[q^{1/2}]]$. In particular, 
  the correction form  $\gamma_{\Theta(-\xi)}$ defined in 
\cite[Def. 4.16]{2011arXiv1103.4217B} vanishes.

\bigskip

Using \cite[Prop. 2.6]{MR2153079} we can define an evaluation
$$\ev_{\Theta(-\xi)}:Q^{\R}_{4m-1}(MSpin)\to  \frac{\R[[q^{1/2}]]}{\Z [[q^{1/2}]]+   \MF^{\Gamma^{0}(2),\R}_{2m} [[q^{1/2}]]}\ .$$
By the intrinsic formula \cite[Thm. 4.17]{2011arXiv1103.4217B}  the evaluation
$$\ev_{\Theta(-\xi)}(\eta^{an}([M,f]))\in \frac{\R[[q^{1/2}]]}{\Z [[q^{1/2}]]+  \MF^{\Gamma^{0}(2),\R}_{2m} [[q^{1/2}]]} $$
 is  represented by the formal power series of reduced $\eta$-invariants
$$-\xi^{\R}(\Dirac\otimes \Theta(\btmr))\in \R[[q^{1/2}]]\ .$$
We now combine  the index theorem \cite[Thm. 3.6]{2011arXiv1103.4217B} stating that $\eta^{an}=\eta^{top}$  and  Lemma  \ref{fwioelfewf}  in order to conclude  that $\eta^{an}([M,f])=0$. Consequently 
 $$\xi^{\R}(\Dirac\otimes \Theta(\btmr ))\in \Z[[q^{1/2}]]+  \MF^{\Gamma^{0}(2),\R}_{2m} [[q^{1/2}]]\ .$$
This is exactly the assertion of  Theorem \ref{thm1000}. \hB

\bibliographystyle{alpha}

\begin{thebibliography}{{Bun}12}

\bibitem[ABP67]{MR0219077}
D.~W. Anderson, E.~H. Brown, Jr., and F.~P. Peterson.
\newblock The structure of the {S}pin cobordism ring.
\newblock {\em Ann. of Math. (2)}, 86:271--298, 1967.

\bibitem[AH68]{MR0231369}
D.~W. Anderson and L.~Hodgkin.
\newblock The {$K$}-theory of {E}ilenberg-{M}ac{L}ane complexes.
\newblock {\em Topology}, 7:317--329, 1968.

\bibitem[APS75]{MR0397797}
M.~F. Atiyah, V.~K. Patodi, and I.~M. Singer.
\newblock Spectral asymmetry and {R}iemannian geometry. {I}.
\newblock {\em Math. Proc. Cambridge Philos. Soc.}, 77:43--69, 1975.

\bibitem[AS69]{MR0259946}
M.~F. Atiyah and G.~B. Segal.
\newblock Equivariant {$K$}-theory and completion.
\newblock {\em J. Differential Geometry}, 3:1--18, 1969.

\bibitem[BD82]{MR679698}
Paul Baum and Ronald~G. Douglas.
\newblock {$K$} homology and index theory.
\newblock In {\em Operator algebras and applications, {P}art {I} ({K}ingston,
  {O}nt., 1980)}, volume~38 of {\em Proc. Sympos. Pure Math.}, pages 117--173.
  Amer. Math. Soc., Providence, R.I., 1982.

\bibitem[{Bec}13]{2013arXiv1310.2851B}
C.~{Becker}.
\newblock {Relative differential cohomology}.
\newblock {\em ArXiv e-prints}, October 2013.

\bibitem[BHS07]{MR2330153}
Paul Baum, Nigel Higson, and Thomas Schick.
\newblock On the equivalence of geometric and analytic {$K$}-homology.
\newblock {\em Pure Appl. Math. Q.}, 3(1, part 3):1--24, 2007.

\bibitem[BM14]{bm}
A.~{Baker} and J.~{Morava}.
\newblock {\$MSp\$ localized away from \$2\$ and odd formal group laws}, March
  2014.

\bibitem[BN10]{MR2652438}
Ulrich Bunke and Niko Naumann.
\newblock The {$f$}-invariant and index theory.
\newblock {\em Manuscripta Math.}, 132(3-4):365--397, 2010.

\bibitem[BNV13]{2013arXiv1311.3188B}
U.~{Bunke}, T.~{Nikolaus}, and M.~{V{\"o}lkl}.
\newblock {Differential cohomology theories as sheaves of spectra}.
\newblock {\em ArXiv e-prints}, November 2013.

\bibitem[Boa95]{board}
J.M. Boardman.
\newblock {\em Handbook of Algebraic Topology}, chapter Stable Operations in
  Generalized Cohomology.
\newblock Elsevier (Amsterdam), 1995.

\bibitem[Bou79]{MR551009}
A.~K. Bousfield.
\newblock The localization of spectra with respect to homology.
\newblock {\em Topology}, 18(4):257--281, 1979.

\bibitem[BS09]{MR2664467}
U.~Bunke and T.~Schick.
\newblock Smooth {$K$}-theory.
\newblock {\em Ast\'erisque}, (328):45--135 (2010), 2009.

\bibitem[BS10]{MR2608479}
U.~Bunke and Th. Schick.
\newblock Uniqueness of smooth extensions of generalized cohomology theories.
\newblock {\em J. Topol.}, 3(1):110--156, 2010.

\bibitem[Bun95]{MR1348799}
U.~Bunke.
\newblock A {$K$}-theoretic relative index theorem and {C}allias-type {D}irac
  operators.
\newblock {\em Math. Ann.}, 303(2):241--279, 1995.

\bibitem[Bun09]{MR2191484}
U.~Bunke.
\newblock Index theory, eta forms, and {D}eligne cohomology.
\newblock {\em Mem. Amer. Math. Soc.}, 198(928):vi+120, 2009.

\bibitem[{Bun}11]{2011arXiv1103.4217B}
U.~{Bunke}.
\newblock {On the topological contents of eta invariants}.
\newblock {\em ArXiv e-prints}, March 2011.

\bibitem[{Bun}12]{2012arXiv1208.3961B}
U.~{Bunke}.
\newblock {Differential cohomology}.
\newblock {\em ArXiv e-prints}, August 2012.

\bibitem[{Fer}14]{2014arXiv1401.1029F}
F.~{Ferrari Ruffino}.
\newblock {Relative (generalized) differential cohomology}.
\newblock {\em ArXiv e-prints}, January 2014.

\bibitem[HS05]{MR2192936}
M.J. Hopkins and I.M. Singer.
\newblock Quadratic functions in geometry, topology, and {M}-theory.
\newblock {\em J. Differential Geom.}, 70(3):329--452, 2005.

\bibitem[HZ04]{MR2153079}
Fei Han and Weiping Zhang.
\newblock Modular invariance, characteristic numbers and {$\eta$} invariants.
\newblock {\em J. Differential Geom.}, 67(2):257--288, 2004.

\bibitem[HZ13]{2013arXiv1312.7494H}
F.~{Han} and W.~{Zhang}.
\newblock { $\eta$-invariant and Modular Forms}.
\newblock {\em ArXiv e-prints}, December 2013.

\bibitem[Lau99]{MR1660325}
Gerd Laures.
\newblock The topological {$q$}-expansion principle.
\newblock {\em Topology}, 38(2):387--425, 1999.

\bibitem[Lau00]{MR1781277}
Gerd Laures.
\newblock On cobordism of manifolds with corners.
\newblock {\em Trans. Amer. Math. Soc.}, 352(12):5667--5688 (electronic), 2000.

\bibitem[LM89]{MR1031992}
H.~B. Lawson, Jr. and M.-L. Michelsohn.
\newblock {\em Spin geometry}, volume~38 of {\em Princeton Mathematical
  Series}.
\newblock Princeton University Press, Princeton, NJ, 1989.

\bibitem[SS10]{MR2732065}
J.~Simons and D.~Sullivan.
\newblock Structured vector bundles define differential {$K$}-theory.
\newblock In {\em Quanta of maths}, volume~11 of {\em Clay Math. Proc.}, pages
  579--599. Amer. Math. Soc., Providence, RI, 2010.

\bibitem[{von}08]{2008arXiv0808.0428V}
H.~{von Bodecker}.
\newblock {On the geometry of the f-invariant}.
\newblock {arXiv0808.0428}, August 2008.

\end{thebibliography}

%
%

\end{document}